\documentclass[leqno]{article}
\usepackage{latexsym,amsmath,amsfonts,mathrsfs,amssymb,amsthm}
\usepackage{yhmath}
\usepackage[usenames]{color}
\usepackage{graphicx}
\usepackage{yhmath}
\def\R{\mathbb R}
\def\N{\mathbb N}

\def\Z{\mathbb Z}

\def\a{\alpha}
\def\e{\epsilon}
\def\d{\delta}
\def\pa{\partial}
\def\Y{\mathbb Y}
\def\T{\mathbb T}

\def\H{{\cal H} }
\def\K{{\cal K} }
\def\mC{\mathfrak C}
\def\mA{\mathfrak A}
\def\mc{\mathfrak c}

\def\mS{\mathfrak S}
\def\ms{\mathfrak s}
\def\cU{{\cal U}}

\def\G{\Gamma}
\def\F{{\cal F} }
\def\oF{\overline{\cal F} }
\def\tb{\textcolor{blue}}

\def\t{\tilde}
\def\lf{\lfloor}
\def\be{\begin{equation}}
\def\ee{\end{equation}}
\def\bs{\backslash}
\def\qed{\hfill$\Box$\bigskip}
\def\nd{\noindent Proof. }
\numberwithin{equation}{section}
\newtheorem{thm}{Theorem}[section]
\newtheorem{lem}[thm]{Lemma}
\newtheorem{pro}[thm]{Proposition}
\newtheorem{defn}[thm]{Definition}
\newtheorem{cor}[thm]{Corollary}
\newtheorem{rem}[thm]{Remark}
\begin{document}
\bigskip

\begin{center}\Large \textbf{Measure and sliding stability for 2-dimensional minimal cones in Euclidean spaces}\end{center}

\bigskip

\centerline{\large Xiangyu Liang}

%

\vskip 1cm

\centerline {\large\textbf{Abstract.}}In this article we prove the measure stability for all 2-dimensional Almgren minimal cones in $\R^n$, and the Almgren (resp. topological) sliding stability for the 2-dimensional Almgren (resp. topological) minimal cones in $\R^3$. As proved in \cite{2T}, when several 2-dimensional Almgren (resp. topological) minimal cones are measure and Almgren (resp. topological) sliding stable, and Almgren (resp. topological) unique, the almost orthogonal union of them stays minimal. As consequence, the results of this article, together with the uniqueness properties proved in \cite{uniquePYT}, permit us to use all 2-dimensional minimal cones in $\R^3$ to generate new families of minimal cones by taking their almost orthogonal unions.

\bigskip


\bigskip


\section{Introduction}

In this article we prove the measure stability for all 2-dimensional Almgren minimal cones in $\R^n$, and the Almgren (resp. topological) sliding stability for all the 2-dimensional Almgren (resp. topological) minimal cones in $\R^3$. The very original motivation of these results comes from the classification of singularities for minimal sets.

The notion of minimal sets (in the sense of Almgren \cite{Al76}, Reifenberg \cite{Rei60}. See David \cite{DJT}, Liang \cite{topo}, etc..for other variances) is a way to try to solve Plateau's problem in the setting of sets. Plateau's problem, as one of the main interests in geometric measure theory, aims at understanding the existence, regularity and local structure of physical objects that minimize the area while spanning a given boundary, such as soap films.

It is known (cf. Almgren \cite{Al76}, David \& Semmes \cite{DS00}) that a $d$-dimensional minimal set $E$ admits a unique tangent plane at almost every point $x$. In this case the local structure around such points are very clear:  the set $E$ is locally a minimal surface (and hence real analytic) around such points, due to the famous result of Allard \cite{All72}. 

So we are mostly interested in what happens around points that admit no tangent plane, namely, the singular points. 

In \cite{DJT}, David proved that the blow-up limits (''tangent objects'') of $d$-dimensional minimal sets at a point are $d$-dimensional minimal cones (minimal sets that are cones in the means time). Blow-up limits of a set a a point reflect the asymptotic behavior of the set at infinitesimal scales around this point. As consequence, a first step to study local structures of minimal sets, is to classify all possible type of singularities--that is to say, minimal cones. 

The plan for the list of $d$-dimensional minimal cones in $\R^n$ is very far from clear. Even for $d=2$, we know very little, except for the case in $\R^3$, where Jean Taylor \cite{Ta} gave a complete classification in 1976, and the list is in fact already known a century ago in other circumstances (see \cite{La} and \cite{He}). They are, modulo isomorphism: a plane, a $\Y$ set (the union of three half planes that meet along a straight line where they make angles of $120^\circ$), and a $\T$ set (the cone over the 1-skeleton of a regular tetrahedron centred at the origin). See the pictures below.

\centerline{\includegraphics[width=0.16\textwidth]{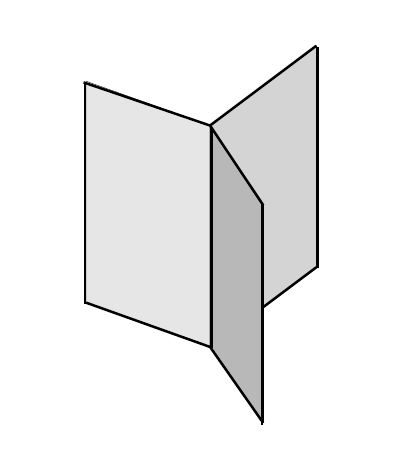} \hskip 2cm\includegraphics[width=0.2\textwidth]{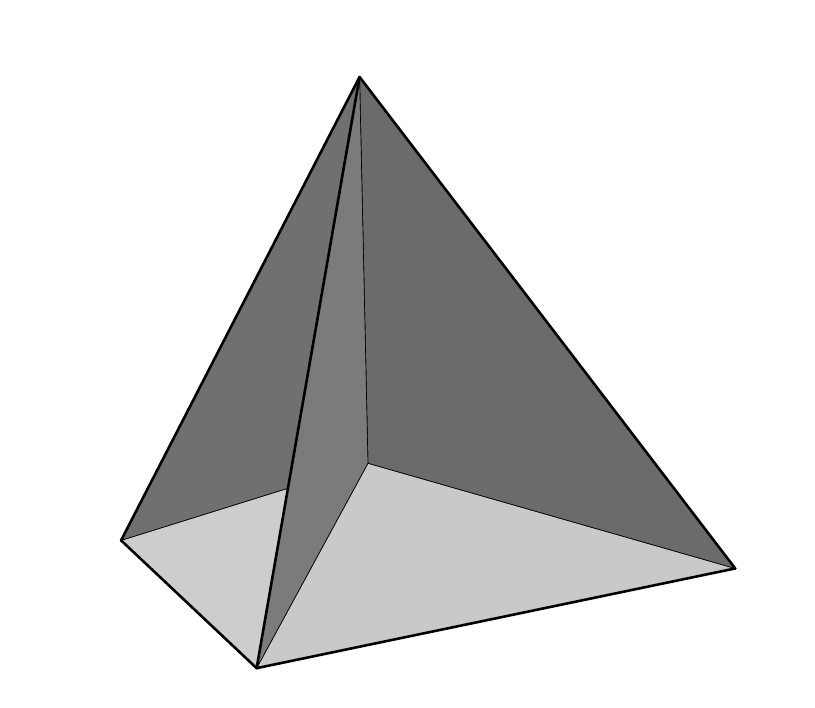}}
\nopagebreak[4]
\hskip 4cm a $\Y$ set\hskip 3.9cm  a $\T$ set

Based on the above, a natural way to find new types of singularities, is by taking unions and products of known minimal cones.

For unions: The minimality of the union of two orthogonal minimal sets of dimension $d$ can be obtained easily from a well known geometric lemma (cf. for example Lemma 5.2 of \cite{Mo84}). Thus one suspects that if the angle between two minimal sets is not far from orthogonal, the union of them might also be minimal.

In case of planes, the author proved in \cite{2p} and \cite{2ptopo}, that the almost orthogonal union of several $d$-dimensional planes is Almgren and topological minimal. When the number of planes is two, this is part of Morgan's conjecture in \cite{Mo93} on the angle condition under which a union of two planes is minimal. 

As for minimal cones other than unions of planes, since they are all with non isolated singularities (after the structure Theorem \tb{2.9}), the situation is much more complicated, as briefly stated in the introduction of \cite{2T}. Up to now we are able to treat a big part of 2 dimensional cases: in \cite{2T} we prove that the almost orthogonal union of several 2-dimensional minimal cones in (in any dimension) are minimal, provided that they are all measure and sliding stable, and satisfy some uniqueness condition. (The theorem is stated separately in the Almgren case and topological case in \cite{2T}.) Moreover, this union stays measure and sliding stable, and satisfies the same uniqueness condition. This enables us to continue obtaining infinitely many new families of minimal cones by taking a finite number of iterations of almost orthogonal unions.

The result makes good sense, because almost all known 2-dimensional minimal cones satisfy all the above conditions. Let us give an account: 

In this article we first prove that all 2-dimensional minimal cones satisfy the measure stability (Theorem \tb{3.1});

We also prove that all 2-dimensional minimal cones in $\R^3$ satisfy the topological and Almgren sliding stability. (Theorem \tb{5.1, 5.2 and 5.11}).

By Theorem \tb{10.1} and Remark \tb{10.5} of \cite{2T}, the almost orthogonal unions of several planes in $\R^n$ are also topological sliding and Almgren sliding stable.

The only known 2-dimensional minimal cone other than the aboves, is the set $Y\times Y$, the product of two 1-dimensional $\Y$ sets. The proof of its sliding minimality is much more involved, so that we will treat it in a separated paper \tb{\cite{stableYXY}}.

We also discuss the uniqueness property for known 2-dimensional minimal cones in \tb{\cite{uniquePYT} and \cite{stableYXY}}.

As a small remark, compare to the unions, the case of product is much more mysterious. It is not known in general whether the product of two non trivial minimal cones stays minimal. We even do not know whether the product of a minimal cone with a line stays minimal. Moreover, if we consider the product of two concrete minimal cones (other than planes) one by one, up to now the only known result is the minimality of the product of two 1-dimensional $\Y$ sets (cf. \cite{YXY}). Among all singular minimal cones, 1-dimensional $\Y$ sets are of simplest structure, but still, the proof of the minimality of their product is surprisingly hard.

Let us say several words for the ingredient of proofs. For the measure stability of any minimal cone, we mainly use second order estimates of the variant of measure for change of center while fixing the boundary, and sliding of boundary while fixing the center; For the sliding properties, we repeatly apply separation and conectness condition to control the measure of sets, and also combine with paired calibration method.

The organization of the rest of the article is the following: 

In Section 2 we introduce basic definitions and preliminaries for minimal sets, and properties for 2-dimensional minimal cones.

We prove the measure stability for 2-dimensional minimal cones in Section 3.

In Section 4 we prove two useful propositions for sliding stabilities, and use them to treat the sliding stabilities for each 2-dimensional minimal cones in $\R^3$ in Section 5.

\textbf{Acknowledgement:} This work is supported by China's Recruitement Program of Global Experts, School of Mathematics and Systems Science, Beihang University. 

\section{Definitions and preliminaries}

\subsection{Some useful notation}

$[a,b]$ is the line segment with end points $a$ and $b$;

$\overrightarrow{ab}$ is the vector $b-a$;

$R_{a,b}$ denote the half line issued from $a$ and passing through $b$;

$B(x,r)$ is the open ball with radius $r$ and centered on $x$;

$\overline B(x,r)$ is the closed ball with radius $r$ and center $x$;

$\H^d$ is the Hausdorff measure of dimension $d$ ;

$d_H(E,F)=\max\{\sup\{d(y,F):y\in E\},\sup\{d(y,E):y\in F\}\}$ is the Hausdorff distance between two sets $E$ and $F$;

$d_{x,r}$ : the relative distance with respect to the ball $B(x,r)$, is defined by
$$ d_{x,r}(E,F)=\frac 1r\max\{\sup\{d(y,F):y\in E\cap B(x,r)\},\sup\{d(y,E):y\in F\cap B(x,r)\}\};$$

For any (affine) subspace $Q$ of $\R^n$, and $x\in Q$, $r>0$, $B_Q(x,r)$ stands for $B(x,r)\cap Q$, the open ball in $Q$;

For any subset $E$ of $\R^n$ and any $r>0$, we call $B(E,r):=\{x\in \R^n: dist (x,E)<r\}$ the $r$ neighborhood of $E$;

For any $d\le n$, any abelien group $G$, and any subset $E\subset \R^n$, $H_d(E,G)$ denote the $d$-th singular homological group of $E$ with coefficient in the group $G$. Let $H_d(E)$ denote $H_d(E,\Z)$ without specification of the coefficient group;

For any $d$-rectifiable subset $E$ of $\R^n$, and for any $x\in E$, let $T_xE$ denote the tangent plane of $E$ at $x$ if it exists.

\subsection{Basic definitions and notations about minimal sets}

In the next definitions, fix integers $0<d<n$. We first give a general definition for minimal sets. Briefly, a minimal set is a closed set which minimizes the Hausdorff measure among a certain class of competitors. Different choices of classes of competitors give different kinds of minimal sets.

\begin{defn}[Minimal sets]Let $0<d<n$ be integers. Let $U\subset \R^n$ be an open set. A relatively closed set $E\subset U$ is said to be minimal of dimension $d$ in $U$ with respect to the competitor class $\mathscr F$ (which contains $E$) if 
\be \H^d(E\cap B)<\infty\mbox{ for every compact ball }B\subset U,\ee
and
\be \H^d(E\bs F)\le \H^d(F\bs E)\ee
for any competitor $F\in\mathscr F$.
\end{defn}

\begin{defn}[Almgren competitor (Al competitor for short)] Let $E$ be relatively closed in an open subset $U$ of $\R^n$. An Almgren competitor for $E$ is an relatively closed set $F\subset U$ that can be written as $F=\varphi_1(E)$, where $\varphi_t:U\to U,t\in [0,1]$ is a family of continuous mappings such that 
\be \varphi_0(x)=x\mbox{ for }x\in U;\ee
\be\mbox{ the mapping }(t,x)\to\varphi_t(x)\mbox{ of }[0,1]\times U\mbox{ to }U\mbox{ is continuous;}\ee
\be\varphi_1\mbox{ is Lipschitz,}\ee
  and if we set $W_t=\{x\in U\ ;\ \varphi_t(x)\ne x\}$ and $\widehat W=\bigcup_{t\in[0.1]}[W_t\cup\varphi_t(W_t)]$,
then
\be \widehat W\mbox{ is relatively compact in }U.\ee
 
Such a $\varphi_1$ is called a deformation in $U$, and $F$ is also called a deformation of $E$ in $U$.
\end{defn}

\begin{defn}[Almgren minimal sets and minimal cones]
Let $0<d<n$ be integers, $U$ be an open set of $\R^n$. An Almgren minimal set $E$ in $U$ is a minimal set defined in Definition \tb{2.1} while taking the competitor class $\mathscr F$ to be the class of all Almgren competitors for $E$.

An Almgren minimal set which is a cone is called a minimal cone.\end{defn}

For future convenience, we also have the following more general definition:

\begin{defn}Let $U\subset \R^n$ be an open set, and let $E\subset \R^n$ be a closed set (not necessarily contained in $U$). We say that $E$ is Almgren minimal in $U$, if $E\cap U$ is minimal in $U$. A closed set $F\subset \R^n$ is called a  deformation of $E$ in $U$, if $F=(E\bs U)\cup \varphi_1(E\cap U)$, where $\varphi_1$ is a deformation in $U$.
\end{defn}

\begin{rem}Since Almgren minimal sets are more often used, we usually omit the word ''Almgren'' and call them minimal sets.
\end{rem}

\begin{defn}[Topological competitors] Let $G$ be an abelian group. Let $E$ be a closed set in an open domain $U$ of $\R^n$. We say that a closed set $F$ is a $G$-topological competitor of dimension $d$ ($d<n$) of $E$ in $U$, if there exists a convex set $B$ such that $\bar B\subset U$ such that

1) $F\bs B=E\bs B$;

2) For all Euclidean $n-d-1$-sphere $S\subset U\bs(B\cup E)$, if $S$ represents a non-zero element in the singular homology group $H_{n-d-1}(U\bs E; G)$, then it is also non-zero in $H_{n-d-1}(U\bs F;G)$.
We also say that $F$ is a $G$-topological competitor of $E$ in $B$.

When $G=\Z$, we usually omit $\Z$, and say directly topological competitor.
\end{defn}

And Definition \tb{2.1} gives the definition of $G$-topological minimizers in a domain $U$ when we take the competitor class to be the class of $G$-topological competitors of $E$.

The simplest example of a $G$-topological minimal set is a $d-$dimensional plane in $\R^n$.  

\begin{pro}[cf.\cite{topo} Proposition 3.7 and Corollary 3.17]  

$1^\circ$ Let $E\subset \R^n$ be closed. Then for any $d<n$, and any convex set $B$, $B'$ such that $\bar B'\subset B^\circ$, every Almgren competitor of $E$ in $B'$ is a $G$-topological competitor of $E$ in $B$ of dimension $d$.

$2^\circ$ All $G$-topological minimal sets are Almgren minimal in $\R^n$.
\end{pro}

\begin{rem}The notion of (Almgren or $G$-topological) minimal sets does not depend much on the ambient dimension. One can easily check that $E\subset U$ is $d-$dimensional Almgren minimal in $U\subset \R^n$ if and only if $E$ is Almgren minimal in $U\times\R^m\subset\R^{m+n}$, for any integer $m$. The case of $G$-topological minimality is proved in \cite{topo} Proposition 3.18.\end{rem}

\subsection{The associated convex domain and stabilities for 2-dimensional minimal cones}

From now on, we are mostly interested in 2-dimensional minimal cones. In this section we will give the definition of the associated convex domains and stabilities for 2-dimensional minimal cones, based on the following structure theorem for 2-dimensional minimal cones.

\begin{thm}[Structure of 2-dimensional minimal cones in $\R^n$, cf. \cite{DJT} Proposition 14.1] Let $K$ be a reduced 2-dimensional minimal cone in $\R^n$, and let $X=K\cap \partial B(0,1)$. Then $X$ is a finite union of great circles and arcs of great circles $C_j,j\in J$. The arcs $C_j$ can only meet at their endpoints, and each endpoint is a common endpoint of exactly three $C_j$, which meet with $120^\circ$ angles (such an endpoint is called a $\Y$ point in $K\cap \pa B$). In addition, the length of each $C_j$ is at least $\eta_0$, where $\eta_0>0$ depends only on the ambient dimension $n$.
\end{thm}

Next we define the convex domain associated to each 2-dimensional minimal cone $K\subset \R^n$. Denote by B the unit ball of $\R^n$. Then by the above theorem, $K\cap\pa B$ is a union of circles $\{s_j,1\le j\le \mu\}$, and arcs of great circles with only $\Y$ type junctions. Let $\eta_0(K)$ denote the minimum of length of these arcs. It is positive, by Theorem \tb{2.9}. 

Denote by $\{a_j, 1\le j\le m\}$ the set of $\Y$ points in $K\cap \pa B$. For any $\eta$, define the $\eta$-convex domain for $K$ 
\be \cU(K,\eta)=\{x\in B: <x, y>< 1-\eta, \forall y\in K\mbox{ and }<x,a_j><1-2\eta, \forall 1\le j\le m\}\subset \R^n. \ee 

From the definition, we see directly that $\cU=\cU(K,\eta)$ is obtained by "cutting off" some small part of the unit ball $B$. More precisely, we first take the unit ball $B$, then just like peelling an apple, we use a knife to peel a thin band (with width about $2\sqrt \eta$) near the net $K\cap \pa B$. Then after this operation, the ball $B$ stays almost the same, except that near the set $K$, the boundary surface will be a thin cylinderical surface. This is the condition ''$<x, y>< 1-\eta, \forall y\in K$''. Next we turn to the singular points $a_j$: they are isolated, so we make one cut at each point, perpendicular to the radial direction, to get a small planar surface near each $a_j$, of diameter about $4\sqrt \eta$. This follows from the condition ''$<x,a_j><1-2\eta, \forall 1\le j\le m $''. 

\medskip

Now for $1\le j,l\le m$, let $\gamma_{jl}$ denote the arc of great circle that connects $a_j$ and $a_l$, if it exists; otherwise set $\gamma_{jl}=\emptyset$. Set $J=\{(j,l): 1\le j,l\le m\mbox{ and }\gamma_{jl}\ne\emptyset\}$, which is exactly the set of pairs $(j,l)$ such that the $\Y$ points $a_j$ and $a_l$ are connected directly by an arc of great circle on $K$. 
 
Denote by $A_j$ the $n-1$-dimensional planar part centered at $(1-2\eta)a_j$ of $\partial \cU$. That is,  
\be A_j=\{x\in \bar B: <x,a_j>=1-2\eta\mbox{ and }<x,y>\le 1-\eta, \forall y\in K\}.\ee
Let $A=\cup_{1\le j\le m}A_j$.

Set 
\be\Gamma_{jl}=\{x\in \bar B, <x,y>=1-\eta\mbox{ for some }y\in \gamma_{jl}\}\bs \mA,\ee
with $\mA$ being the cone over $A$ centered at 0,
and
\be S_j=\{x\in \bar B,<x,y>=1-\eta\mbox{ for some }y\in s_j\}.\ee
Then $\Gamma_{jl}$ is the band like part of $\pa \cU$ near each $(1-\eta)\gamma_{jl}$, and similar for $S_j$. The union
$\Gamma=\cup_{1\le j,l\le m}\Gamma_{jl}$ together with $S=\cup_{1\le j\le \mu}S_j$ is the whole cylinderical part of $\pa \cU$.

Set $\mathfrak c_{jl}$, $\mC_{jl}$, $\mC$, $\ms_j$, $\mS_j$, $\mS$, $\mA_j$ and $\mA$ the part of the cone (centered at 0) included in $\cU$ over $\gamma_{jl}$, $\Gamma_{jl}$, $\Gamma$, $s_j$, $S_j$, $S$, $A_j$ and $A$ respectively, where for any set $S\subset \R^n$, the cone over $S$ is defined to be $\{ts: x\in S, t\ge 0\}$.

\medskip

Let $\eta_1(K)$ be the superium of the number $\eta$, such that on $\pa\cU(K,\eta)$, any 3 of the $\Gamma_{jl}, (j,l)\in J$ and $S_j,1\le j\le \mu$ never have a common point, and the $A_j,1\le j\le m$ are disjoint.

\medskip

\textbf{In the rest of this article, when we treat each 2-dimensional minimal cones $K$, we will only consider $\eta<\eta_1(K)$.}

\medskip

Under this condition, the shape of the planar region $A_i$ (of dimension $n-1$) will be obained by cutting off 3 small planar part of a $n-1$-dimensional ball . Take $a_i$ for example, suppose, without loss of generality, that the three $\Y$ points in $K\cap \pa B$ that are adjacent to $a_1$ are $a_l,l=2,3,4$. Then the planar region centered at $(1-2\eta)a_1$ is obtained by:

Firstly, take the $n-1$-dimensional ball $\Omega_1$ perpendicular to $\overrightarrow{oa_1}$ and centered at $(1-2\eta)a_1$ (hence the radius of $\Omega_1$ is $R=\sqrt{1-(1-2\eta)^2}$). For $l=2,3,4$, denote by $x_l$ the intersection of $(1-\eta)\gamma_{1l}$. Then the $x_l,l=2,3,4$ will situated on a 2-plane passing through the center of the ball $\Omega_1$, hence they belong to a same great circle. Let $L_{1l}$ be the $n-1$ subspace containing $x_l$ and orthogonal to $\overrightarrow{(1-2\eta)a_1, x_l}$.  We then cut off the small part of  $\Omega_1$ that is on the other side of $L_{1l}, l=2,3,4$, and get the planar region $A_1$ . This forms part of the boundary of $\cU$.

Let us look at the boundary of $A_1$: By definition, it is composed of three disjoint $n-2$-dimensional small balls : 
$$I_{1l}=\{x\in B: <x,a_1>=1-2\eta\mbox{ and }<x,y>=1-\eta\mbox{ for some }y\in \gamma_{1l}\},2\le l\le 4$$
and the rest of the boundary of $\Omega_1$. Note that the diameter of $I_{1l}$ is $R_1=\sqrt{1-(1-\eta)^2}$.

See the figure below. 

\centerline{\includegraphics[width=0.7\textwidth]{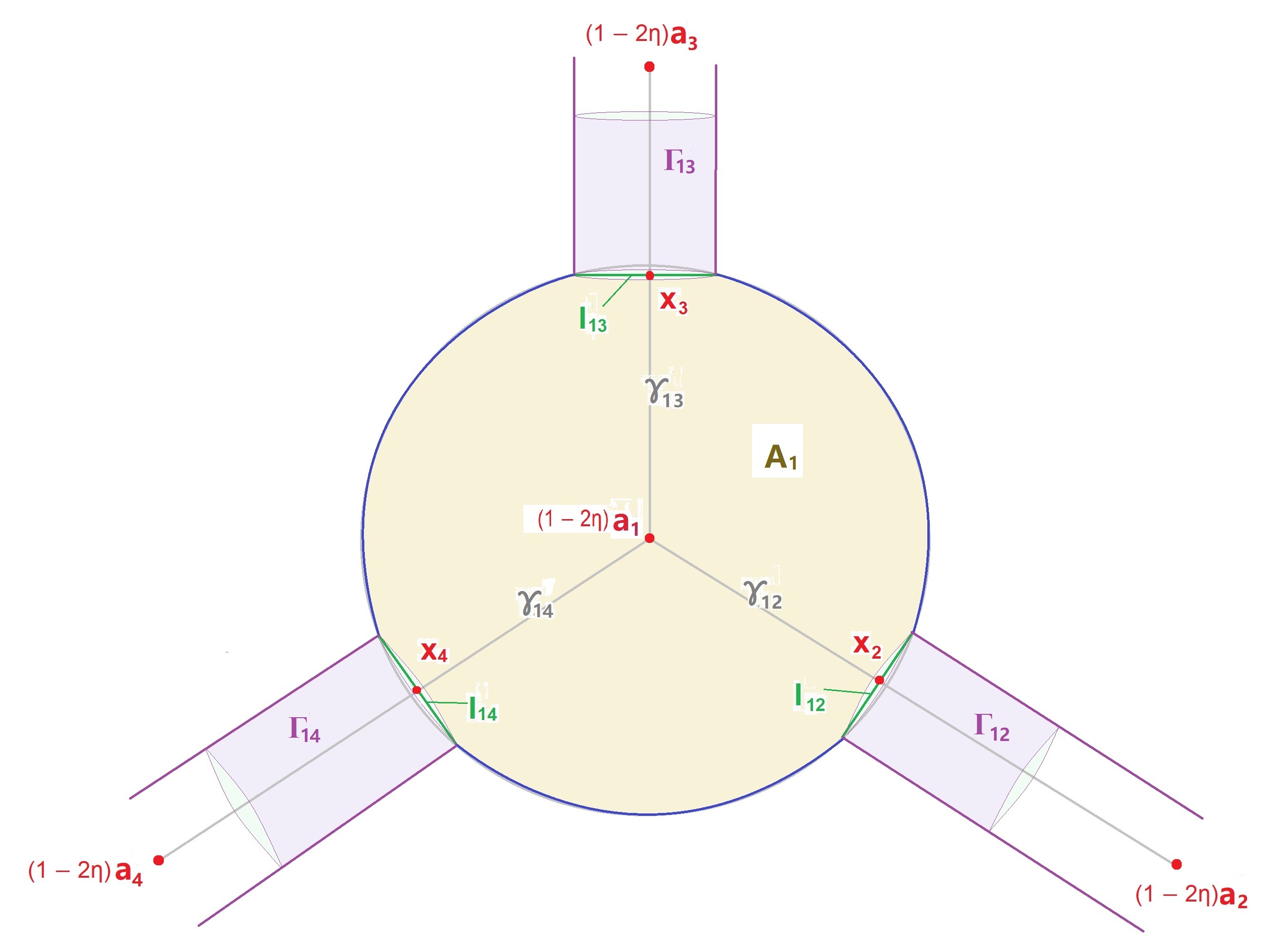} }

We define similarly, for any $(j,l)\in J$: $I_{jl}$. Then the boundary of $A_j, j\in m$ is the union of the balls $I_{jl},(j,l)\in J$, and the rest of the sphere $\pa\Omega_1$.

Set $A=\cup_{1\le j\le m}A_j$. This is the whole planar part of $\pa \cU$. 

For the rest of $\pa \cU$ that is not spherical, they are obtained by the equation
\be x\in B, <x,y>\le 1-\eta, \forall y\in K.\ee

\begin{defn}Let $U$ be an open subset of $\R^n$, let $E\subset \bar U$ be closed. For $\d>0$, a $\d$-sliding Lipschitz deformation of $E$ in $\bar U$ is a set $F\subset \bar U$ that can be written as $F=\varphi_1(E)$, where $\varphi_t: \bar U\to \bar U$ is a family of continuous mappings such that
\be \varphi_0(x)=x\mbox{ for }x\in \bar U;\ee
\be \mbox{ the mapping }(t,x)\to \varphi_t(x)\mbox{ of }[0,1]\times \bar U\mbox{ to }\bar U\mbox{ is continuous};\ee
\be \varphi_1\mbox{ is Lipschitz },\ee
\be \varphi_t(\partial U)\subset\partial U,\ee
and
\be|\varphi_t(x)-x|<\d\mbox{ for all }x\in E\cap\partial U.\ee

Such a $\varphi_1$ is called a sliding deformation in $\bar U$, and $F$ is called a $\d$-sliding deformation of $E$ in $\bar U$.
\end{defn}

Let $\F_\d(E,\bar U)$ denote the set of all $\d$-sliding deformation of $E$ in $\bar U$, and let $\oF_\d(E,\bar U)$ be the family of sets that are Hausdorff limits of sequences in $\F_\d(E,\bar U)$. That is: we set
\be \begin{split}\oF_\d(E,U)&=\{F\subset\bar U: \tb{(2.1)}\mbox{ holds for }F\mbox{, and }\exists \{E_n\}_n\subset\F_\d(E,U)\mbox{ such that }d_H(E_n,F)\to 0
\}.\end{split}\ee

\begin{defn}[$\d$-(Almgren) sliding minimal sets]Let $\d>0$, $U\subset \R^n$ be open, and let $E\subset \bar U$ be closed. We say that $E$ is $\d$-(Almgren) sliding minimal in $\bar U$, if \tb{(2.1)} holds, and \tb{(2.2)} holds for all $F\in \oF_\d(E,\bar U)$.
\end{defn}

\begin{defn}[Stable minimal cones] Let $K$ be a 2-dimensional Almgren minimal cone in $\R^n$.

$1^\circ$ We say that $K$ is $(\eta, \d)$-Almgren sliding stable, if for some $\eta\in (0,\eta_1(K))$, and $\d\in (0,\eta)$, $K$ is $\d$-Almgren sliding minimal in $\bar \cU(K,\eta)$. We say that $K$ is Almgren sliding stable if it is $(\eta, \d)$-sliding stable for some $\eta\in (0,\eta_1(K))$ and $\d\in (0,\eta)$.

$2^\circ$ We say that $K$ is $(\eta, \d)$-measure stable, if there exists $\eta\in (0,\eta_1(K))$, and $\d\in(0,\eta)$, such that for all $y\in \R^n$ with $||y||<\d$, we have
\be \H^2(K\cap \cU(K,\eta))=\H^2((K+y)\cap \cU(K,\eta)).\ee
We say that $K$ is measure stable if it is $(\eta, \d)$-measure stable for some $\eta\in (0,\eta_1(K))$ and $\d\in(0,\eta)$.

\end{defn}

\begin{defn}Let $K$ be a 2-dimensional $G$-topological minimal cone in $\R^n$. Let $0<\d<\eta<\eta_1(K)$. 
$1^\circ$ We say that a closed set $F$ is an $(\eta,\d)$-$G$-topological sliding competitor for $K$, if there exists a 2-dimensional $G$-topological competitor $E$ of $K$ in $\cU(K,\eta)$, such that $F$ is a $\d$-sliding deformation of $E$ in $\bar\cU(K,\eta)$.

$2^\circ$ We say that $K$ is $(\eta,\d)$-$G$-topological sliding stable, if for all $(\eta,\d)$-$G$-topological sliding competitor $F$ of $K$, we have
\be \H^d(F\cap \bar\cU(K,\eta))\ge \H^d(K\cap \bar\cU(K,\eta)).\ee
\end{defn}

For the proof of Almgren and $G$-topological sliding stability, we also introduce the following stronger stability property:

\begin{defn}Let $K$ be a 2-dimensional $G$-topological minimal cone in $\R^n$. Let $0<\d<\eta<\eta_1(K)$, and let 
\be V_\d:=\{x\in \partial \cU(K,\eta): \mbox{dist}(x,K)\le \d\}.\ee

$1^\circ$ We say that a closed set $F$ is an $(\eta,\d)$ $G$-topological competitor for $K$, if 
\be F\bs \bar\cU(K,\eta)=K\bs \bar\cU(K,\eta),\ee 
\be F\cap \pa\cU(K,\eta)\subset V_\d,\ee
and
\be F\cup V_\d\mbox{ is a }G\mbox{-topological competitor for }K.\ee 

$2^\circ$ We say that $K$ is $(\eta,\d)$-$G$-topological stable, if for all $(\eta,\d)$ $G$-topological competitor $F$ of $K$, we have
\be \H^d(F\cap\bar\cU(K,\eta))\ge \H^d(K\cap \bar\cU(K,\eta)).\ee
\end{defn}

\begin{rem}$1^\circ$ It is easy to see that, for any $\eta<\eta'$ and $\d<\d'$, $K$ is $(\eta',\d')$-Almgren (resp. $G$-topological) sliding stable $\Rightarrow K$ is $(\eta,\d)$-Almgren (resp. $G$-topological) sliding stable. Same for the $G$-topological stability.

$2^\circ$ In contrast to the definition of sliding stabilities, in the definition of the measure stability we see no difference between Almgren and topological minimality.

$3^\circ$  Let $K^1$ and $K^2$ be two $(\eta,\d)$-measure stable minimal cones of dimension 2, with dist$(K^1\cap \pa B,K^2\cap \pa B)$ relatively large, then even if the cone $C=K^1\cup K^2$ might not be minimal, we can define $\cU(C,\eta)$, and the measure stability in the sense definition \tb{2.12}. Directly from the definition, we know that $C$ is also $(\eta,\d)$-measure stable.
 \end{rem}

\section{Measure stability of 2-dimensional minimal cones}

In this section we prove the measure stability for all 2-dimensional minimal cones in $\R^n$.

\begin{thm}[Measure stability for 2-dimensional minimal cones]
Let $K$ be a 2-dimensional minimal cone in $\R^n$, then for each $\eta<\eta_1(K))$, $K$ is $(\eta, \eta)$-measure stable. That is, let $\cU=\cU(K,\eta)$ be the $\eta$-convex domain associated to $K$, then for all $q\in \R^n$ with $||q||<\eta$, we have
\be \H^2(K\cap \cU)=\H^2((K+q)\cap \cU).\ee
\end{thm}

\nd It is enough to prove that, for any fixed $\eta<\eta_1(K)$, any $q\in\partial B(0,1)$, and any $|t_0|<\eta$,
\be \frac{d}{dt}|_{t=t_0}\H^2((K+tq)\cap \cU(K,\eta))=0.\ee

So fix any $\eta$ and $q$ as above. Let $\cU$ denote $\cU(K,\eta)$. For $t<\eta$, let $K_t$ denote $(K+tq)\cap \cU$, and set $X_t=K_t\cap \partial \cU$. Then $K_0=K\cap\bar\cU$.

Fix any $t_0$ with $|t_0|<\eta$. Fix $s$ small so that $|t_0\pm s|<\eta$. Set $o_s:=(t_0+s)q$. For $t$ with $|t_0\pm t|<\eta$, let $C_t\subset \cU$ be the cone over $X_{t_0+t}$ centered at $o_s$. Then $C_s=K_{t_0+s}$.

Let us first estimate the difference of measure between $C_0$ and $C_s=K_{t_0+s}$. They are cones with the same center $o_s$. 

Let $\mathfrak c^s_{jl}$, $\mC^s_{jl}$, $\mC^ss$, $\ms^s_j$, $\mS^s_j$, $\mS^s$, $\mA^s_j$ and $\mA^s$ denote the part of the cone centered at $o_s$ included in $\cU$ over $\gamma_{jl}$, $\Gamma_{jl}$, $\Gamma$, $s_j$, $S_j$, $S$, $A_j$ and $A$ respectively.

By definition, we know that on $\pa\cU$,  for any $|t|<\eta$, $B(X_t, \eta)\subset A\cup\G\cup S$. Therefore, since $|t_0|<\eta$, for $t$ small, we have
\be X_{t_0+t}\subset  B(X_0, \eta)\subset A\cup\G\cup S=[\cup_{1\le i\le m}A]\cup[\cup_{(j,l)\in J}\G_{jl}]\cup[\cup_{1\le i\le \mu} S_i],\ee
and thus, 
\be C_t\subset \mA^s\cup \mC^s\cup \mS^s=[\cup_{1\le i\le m}\mA_i^s]\cup[\cup_{(j,l)\in J}\mC_{jl}^s]\cup[\cup_{1\le i\le \mu}\mS_i^s],\ee
where the unions on the right-hand-sides of the above two equations are both disjoint. 

$1^\circ$ In the regions $\mA_i^s, 1\le i\le m$. Take $i=1$ for example. 
For $|t|<\eta$, et $Y_t$ denote $A_1\cap X_t$, it is a 1-dimensional $\Y$ set in $A_1$. 

Let $I_1, I_2,I_3$ denote the three $n-1$-dimensional disks $I_{1j}, (1,j)\in J$ for short, and let $y_\a$ denote the the centers of $I_\a$, $1\le \a\le 3$. Then $y_\a, 1\le \a\le 3$ are also the three points of intersection of $Y_0$ with $\pa A_1$. The center of $A_1$ is $(1-2\eta)a_1$.

Let $P$ be the 2-dimensional affine subspace containing $Y_0$, and let $\pi$ denote the orthogonal projection from $\R^n$ to $P$. Let $b$ denote $\pi((1-\eta)a_1)=\pi(a_1)$.

The projection of $A_1$ is the 2-dimensional convex region
\be \pi(A_1)=B(b, R)\bs (\cup_{\a=1}^3\{x\in B(b,R)\cap P:<x, a_\a>>1-\eta\}).\ee

For $\a=1,2,3$, let $S_\a$ denote $\pi(I_\a)$, then it is a segment centered at $a_\a$ with length $R_1=2\sqrt{1-(1-\eta)^2}$. The boundary of $\pi(A_1)$ is a union of three segments $S_\a, 1\le \a\le 3$, and three arcs of circles $\xi_\a\subset \pa B(b,R)$, $\a=1,2,3$. Let $L_\a$ denote the line containing $S_\a, 1\le \a\le 3$. See the picture below. Then the angles between any two of the $L_\a$ is $\frac \pi 3$, and they enclose a equilateral triangle $\Delta$, whose center is $b$.

\centerline{\includegraphics[width=0.7\textwidth]{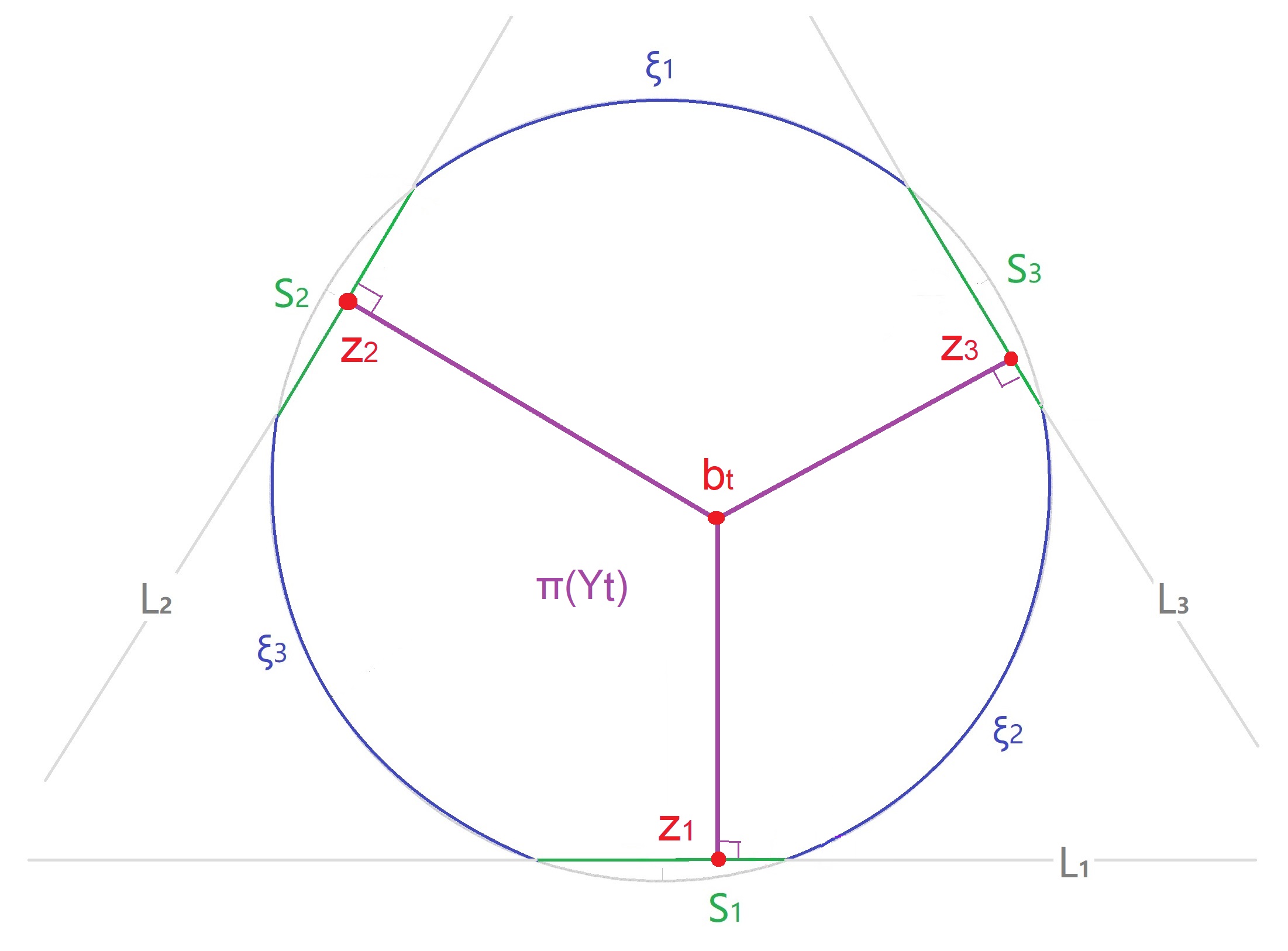} }

Now for each $|t|<\eta$, we know that $Y_t$ is parallel to $Y_0$, hence $\pi(Y_t)$ is also a 1-dimensional $\Y$ set, and
\be \H^1(\pi(Y_t))=\H^1(Y_t).\ee
 On the other hand, since $Y_t\subset B(Y_0, \eta)$, the intersection $Y_t\cap \pa \pi(A_1)$ is contained in $\cup_{1\le \a\le 3}S_\a$, and $\pi(Y_t)$ intersects each $S_\a$ precisely at one point $z_\a$. Let $b_t$ denote the center of $\pi(Y_t)$. Then the segment $[b_t, z_\a]$ is perpendicular to $S_\a$, and hence $L_\a$. As a result, we know that 
\be \H^1(\pi(Y_t))=\sum_{1\le\a\le 3}|z_\a-b_t|=\sum_{1\le\a\le 3}\mbox{dist}(b_t, L_\a).\ee

On the other hand, since $\Delta$ is a equilateral triangle, we know that for any point $p\in \Delta$, the quantity $\sum_{1\le\a\le 3}\mbox{dist}(p, L_\a)$ is constant. As result, since $b_t\subset B(b, \eta)\subset \Delta^\circ$, we know that 
\be \sum_{1\le\a\le 3}\mbox{dist}(b_t, L_\a)=\sum_{1\le\a\le 3}\mbox{dist}(b, L_\a)=\H^1(Y_0).\ee

Combine \tb{(3.6), (3.7) and (3.8)}, we get
\be \H^1(Y_t)=\H^1(Y_0)\ee
for all $-\eta<t<\eta$.

Now let us look at the measure of $C_s\cap \mA_1^s$ and $C_0\cap \mA_1^s$. They are cones centered at the same point $o_s$, that is, 
\be C_s\cap \mA_1^s=(Y_{t_0+s}\# o_s)\cap \mA_1^s\mbox{ and }C_0\cap \mA_1^s=(Y_{t_0}\# o_s)\cap \mA_1^s.\ee

Set $f:\R^n\to \R: f(x)=<x,a_1>$, then $P_1=ker f$ is a $n-1$ dimensional subspace orthogonal to the 1-subspace $L$ generated by $\vec a_1$. Write $\R^n=P_1\times L$. Let $\pi_1$ denote the orthogonal projection from $\R^n$ to $P_1$, and set $Y=\pi_1(Y_{t_0})$. Let $g$ be the restriction of $f$ on $C_0\cap \mA_1^s$. By the structure of $Y_t$, which is smooth except for the only singular point $\xi_t$, we know that for all $x\in C_0\cap A_1^s\bs [o_s,\xi_{t_0}]$, the tangent plane $T_x$ of $C_0\cap \mA_1^s$ at $x$ exists. For a such point $x$, the differential $Dg(x)$ of $g$ at $x$ is a linear map from the 2-dimensional subspace $T_x$ to $\R$.
Hence by Coarea formula (\cite{Fe} 3.2.22), we know that
\be \begin{split}&\H^2(C_0\cap \mA_1^s)=\int_{f(o_s)}^{1-2\eta}dr \int_{g^{-1}\{r\}\cap C_0\cap \mA_1^s}||D_g(x)||^{-1}d\H^1(x)\\
&=\int_{f(o_s)}^{1-2\eta}dr \int_{[(\frac{r-f(o_s)}{1-2\eta-f(o_s)})(Y-\pi_1(o_s))+\pi_1(o_s)]\times\{r a_1\}}||D_g(x)||^{-1}d\H^1(x)\\
&=\int_{f(o_s)}^{1-2\eta}dr (\frac{r-f(o_s)}{1-2\eta-f(o_s)})\int_Y||D_g(([(\frac{r-f(o_s)}{1-2\eta-f(o_s)})(y-\pi_1(o_s))+\pi_1(o_s)], r))||^{-1}d\H^1(y)
\end{split}\ee
 Now fix any $y\in Y\bs \{\pi_1(\xi_{t_0})\}$, and any $r\in [f(o_s),1-2\eta]$, since $C_0$ is the cone over $Y_{t_0}$ centered at $o_s$, we know that the tangent plane of $C_0$ at the point $([(\frac{r-f(o_s)}{1-2\eta-f(o_s)})(y-\pi_1(o_s))+\pi_1(o_s)], r)$ does not depend on $r$, and is equal to the 2-plane $Q_{y,t}$ spaned by $y-\pi_1(\xi_{t_0})$ and $\xi_{t_0}-o_s$.  Therefore the quantity $\int_Y||D_g(([(\frac{r-f(o_s)}{1-2\eta-f(o_s)})(y-\pi_1(o_s))+\pi_1(o_s)], r))||^{-1}d\H^1(y)$ is also independent of $r$, and if we let $\theta_{y,t}$ denote the angle between $Q_{y,t}$ and the line $L$ generated by $a_1$, then 
 \be \begin{split}
||D_g(([(\frac{r-f(o_s)}{1-2\eta-f(o_s)})(y-\pi_1(o_s))+\pi_1(o_s)], r))||=\cos\theta_{y,t}\\
 =\max\{<x,a_1>: x\in Q_{y,t}\mbox{ and }||x||=1\}\ge <\frac{\xi_{t_0}-o_s}{||\xi_{t_0}-o_s||},a_1>.\end{split}\ee

 But we know that $\xi_{t_0}\in K_0+t_0q$ is a singular point of $K_{t_0}$, hence $\xi_{t_0}+sq\in K_{t_0+s}$ is also a singular point of $K_{t_0+s}$, and thus $\xi_{t_0}+sq-o_s\in K_0$ is a singular point $z$ of $K_0$. But $\xi_{t_0}+sq-o_s\in \mA_1$, hence should be a multiple of $a_1$. Also note that $\xi_{t_0}-o_s=(\xi_{t_0}-\xi_{t_0+s})+(\xi_{t_0+s}-o_s)$, where $\xi_{t_0}-\xi_{t_0+s}\in P_1$, which is orthogonal to $a_1$, and $\xi_{t_0+s}-o_s$ is a multiple of $a_1$. As a result, 
 \be ||\xi_{t_0}-\xi_{t_0+s}||=\mbox{dist}(\xi_{t_0}-o_s,L)\le ||\xi_{t_0}-o_s-\lambda a_1||, \forall \lambda\in \R.\ee
 
 In particular, 
 \be ||\xi_{t_0}-\xi_{t_0+s}||\le ||\xi_{t_0}-o_s-(\xi_{t_0}+sq-o_s)||=s,\ee
because $\xi_{t_0}+sq-o_s$ is a multiple of $a_1$.

Again, because $\xi_{t_0}-\xi_{t_0+s}$ is perpendicular to $\xi_{t_0+s}-o_s$, we have
\be ||\xi_{t_0}-o_s||^2=||\xi_{t_0}-\xi_{t_0+s}||^2+||\xi_{t_0+s}-o_s||^2\le s^2+||\xi_{t_0+s}-o_s||^2,\ee
and hence
\be \begin{split}||\xi_{t_0}-o_s||&\le\sqrt{s^2+||\xi_{t_0+s}-o_s||^2}=||\xi_{t_0+s}-o_s||\sqrt{1+(\frac{s}{||\xi_{t_0+s}-o_s||})^2}\\
&\le ||\xi_{t_0+s}-o_s||(1+\frac12(\frac{s}{||\xi_{t_0+s}-o_s||})^2)\le ||\xi_{t_0+s}-o_s||+s^2,\end{split}\ee
as $\eta$, $t_0$ and $s$ are all small.

Combine with \tb{(3.12)}, we have
\be \begin{split}
||D_g(([(\frac{r-f(o_s)}{1-2\eta-f(o_s)})&(y-\pi_1(o_s))+\pi_1(o_s)], r))||\ge \frac{1}{||\xi_{t_0}-o_s||}<\xi_{t_0}-o_s, a_1>\\
&=\frac{1}{||\xi_{t_0}-o_s||}||\xi_{t_0+s}-o_s||\ge \frac{||\xi_{t_0+s}-o_s||}{||\xi_{t_0+s}-o_s||+s^2},\end{split}\ee
for $s$ small.

On the other hand, since $g$ is 1-Lipschitz, 
\be ||D_g(([(\frac{r-f(o_s)}{1-2\eta-f(o_s)})(y-\pi_1(o_s))+\pi_1(o_s)], r))||\le 1.\ee

Recall that $([(\frac{r-f(o_s)}{1-2\eta-f(o_s)})(y-\pi_1(o_s))+\pi_1(o_s)], r)$ represents all points $x\in C_0\cap \mA_1^s$ in the coordinate $P_1\times L$, hence 
\be  \frac{||\xi_{t_0+s}-o_s||}{||\xi_{t_0+s}-o_s||+s^2}\le ||Dg(x)||\le 1, \forall x\in  C_t\cap \mA_1.\ee

Now we bring \tb{(3.19)} back to \tb{(3.11)}, and get
\be\begin{split}\int_{f(o_s)}^{1-2\eta}dr &(\frac{r-f(o_s)}{1-2\eta-f(o_s)})\int_Yd\H^1(y)\le \H^2(C_0\cap \mA_1^s)\\
&\le \int_{f(o_s)}^{1-2\eta}dr (\frac{r-f(o_s)}{1-2\eta-f(o_s)})\int_Y\frac{||\xi_{t_0+s}-o_s||+s^2}{||\xi_{t_0+s}-o_s||}d\H^1(y).\end{split}\ee

Therefore
\be \begin{split}0&\le \H^2(C_0\cap \mA_1^s)-\int_{f(o_s)}^{1-2\eta}dr (\frac{r-f(o_s)}{1-2\eta-f(o_s)})\int_Yd\H^1(y)\\
&\le \int_{f(o_s)}^{1-2\eta}dr (\frac{r-f(o_s)}{1-2\eta-f(o_s)})\int_Y(\frac{||\xi_{t_0+s}-o_s||+s^2}{||\xi_{t_0+s}-o_s||}-1)d\H^1(y)\\
& =\int_{f(o_s)}^{1-2\eta}dr (\frac{r-f(o_s)}{1-2\eta-f(o_s)})\int_Y(\frac{s^2}{||\xi_{t_0+s}-o_s||})d\H^1(y).
\end{split}\ee
Note that $||\xi_{t_0+s}-o_s||=1-2\eta-f(o_s)$, hence
\be \begin{split}&\int_{f(o_s)}^{1-2\eta}dr (\frac{r-f(o_s)}{1-2\eta-f(o_s)})\int_Y(\frac{s^2}{||\xi_{t_0+s}-o_s||})d\H^1(y)\\
&=\frac{s^2}{(1-2\eta-f(o_s))^2}\H^1(Y)\int_{f(o_s)}^{1-2\eta}dr(r-f(o_s))\\
&=\frac{s^2}{(1-2\eta-f(o_s))^2}\H^1(Y)[\frac12 (1-2\eta-f(o_s)]^2\\
&=\frac{s^2}{2}\H^1(Y).
\end{split}\ee
 Therefore
\be 0\le \H^2(C_0\cap \mA_1^s)-\int_{f(o_s)}^{1-2\eta}dr \int_Y(\frac{r-f(o_s)}{1-2\eta-f(o_s)})d\H^1(y)\le \frac{s^2}{2}\H^1(Y).\ee

On the other hand, we know that $\H^1(Y)=\H^1(\pi_1(Y_t))=\H^1(Y_t)=\H^1(Y_0)$, and the spine of $C_s$ is parallel to $L$, hence by the same argument above, 
\be \H^2(C_s\cap \mA_1^s)=\int_{f(o_s)}^{1-2\eta}dr \int_Y(\frac{r-f(o_s)}{1-2\eta-f(o_s)})d\H^1(y).\ee
Combine with \tb{(3.23)}, we have
\be 0\le \H^2(C_0\cap \mA_1^s)-\H^2(C_s\cap \mA_1^s)\le \frac{s^2}{2}\H^1(Y_0).\ee

Similarly, for all $1\le i\le m$, we have
\be |\H^2(C_0\cap \mA_i^s)-\H^2(C_s\cap \mA_i^s)|\le \frac{s^2}{2}\H^1(Y_0).\ee
We sum over all $i$, and get
\be |\H^2(C_0\cap \mA^s)-\H^2(C_s\cap \mA^s)|=\frac{m\H^1(Y_0)}{2}s^2.\ee

$2^\circ$ In the regions $\mC^s_{jl}, (j,l)\in J$. Suppose for example $(j,l)=(1,2)$. Let $Q$ denote the $n-2$ dimensional subspace orthogonal to $\mc_{12}$, then $\Gamma_{12}=((1-\eta)\gamma_{12}\bs \mA)\times B_Q(0, R_1)$, where 
$R_1=\sqrt{1-(1-\eta)^2}$ as before.

For $|t|<R_1$, set $X_t=K_t\cap \G_{12}$. Then $X_0$ is just the arc of circle $((1-\eta)\gamma_{12}\bs \mA)$, and $X_t$ will be a translation of $X_0$, and 
\be \mbox{dist}(X_0,X_t)\le |tq|=t.\ee

But since $\G_{12}$ is the product of $X_0$ with the ball $B_Q(0, R_1)$ of its orthogonal space, we know that in $\Gamma_{12}$, all translations within distance $R_1$ of $X_0$ has the same $\H^1$ measure:
\be \H^1(X_0)=\H^1(X_t)\mbox{ for }t<R_1.\ee

Now as before, for $t$ small, we know that $C_t\cap \mC^s_{12}$ is the cone over $X_{t_0+t}$ centered at $o_s$. Since $X_{t_0+t}$ and $X_{t_0+s}$ are both translations of $X_0$ in $\Gamma_{12}=X_0\times B_Q(0, R_1)$, there exists $y_t\in B_Q(0, R_1)$ such that $X_{t_0+t}=X_{t_0+s}+y$. Then 
\be |y_t|=\mbox{dist}(X_{t_0+t},X_{t_0+s})\le |t-s|.\ee

For each $z\in X_{t_0+t}$, $z-y\in X_{t_0+s}$, by definition, we know that $||z-y-o_s||$ is constant, which is equal to $l_s:=1-\eta-<(t_0+s)q, a_1>$. Still since $z-y\in X_{t_0+s}$, $z-y-o_s$ is orthogonal to $y$, hence
\be ||z-o_s||=\sqrt{||y_t||^2+||z-y-o_s||^2}=\sqrt{||y_t||^2+l_s^2}.\ee

As a result, since $C_t\cap \mC^s_{12}$ is a cone over $X_{t_0+t}$ centered at $o_s$ in $\cU$, we know that
\be\H^2(C_t\cap \mC^s_{12})=\frac12\sqrt{||y_t||^2+l_s^2}\H^1(X_{t_0+t})=\frac12\sqrt{||y_t||^2+l_s^2}\H^1(X_0),\ee
by \tb{(3.29)}.

When $t=s$, we know that $|y_s|\le |s-s|=0$, hence
\be \H^2(C_s\cap \mC^s_{12})=\frac12l_s^2\H^1(X_0),\ee
and $y_t\ge y_s$ for $t$ small. Therefore
\be 0\le \H^2(C_0\cap \mC^s_{12})-\H^2(C_s\cap \mC^s_{12})=\frac12(\sqrt{||y_0||^2+l_s^2}-l_s)\H^1(X_0).\ee

Note that $|l_s|=|1-\eta-(t_0+s)<q,a_1>|\ge 1-\eta-|t_0+s|\ge 1-3\eta>\frac12$ for $s<\eta$. Hence
\be \sqrt{||y_0||^2+l_s^2}-l_s=l_s(\sqrt{1+(\frac{||y_0||}{l_s})^2}-1)\le l_s(1+\frac12(\frac{||y_0||}{l_s})^2-1)=\frac{||y_0||^2}{2l_s}\le\frac{s^2}{2l_s}\le s^2.\ee

Therefore, 
\be |\H^2(C_0\cap \mC^s_{12})-\H^2(C_s\cap \mC^s_{12})|\le \frac12s^2\H^1(X_0).\ee

Similar argument gives
\be |\H^2(C_0\cap \mC^s_{jl})-\H^2(C_s\cap \mC^s_{jl})|\le \frac12s^2\H^1(X_0), \forall (j,l)\in J.\ee

We sum over all $(j,l)\in J$, and get
\be |\H^2(C_0\cap \mC^s)-\H^2(C_s\cap \mC^s)|\le \frac{|J|\H^1(X_0)}{2}s^2.\ee
 
 $3^\circ$ In the regions $\mS_j,1\le j\le \mu$. The argument is exactly the same as $2^\circ$. Hence we also have
 \be |\H^2(C_0\cap \mS^s)-\H^2(C_s\cap \mS^s)|\le Cs^2,\ee
 where $C$ is a constant that does not depend on $s$.
 
Now by $1^\circ$-$3^\circ$, and \tb{(3.4)}, we have
\be \H^2(C_0)-\H^2(C_s)=O(s^2).\ee

On the other hand, denote by $G_s$ the cone over $X_{t_0}$ centered at $o_s=(t_0+s)q$ and contained in $\cU$. Then by minimality of $K_{t_0}$, we know that 
\be\frac{d}{ds}|_{s=0}\H^2(G_s)=0,\ee
and hence 
\be |\H^2(G_s)-\H^2(K_{t_0})|=O(s^2).\ee

Note that $C_s=K_{t_0+s}$ and $C_0=G_s$, hence \tb{(3.40) and (3.42)} yields
\be |\H^2(K_{t_0+s})-\H^2(K_{t_0})|=O(s^2).\ee

As a result, for $t_0<\eta$, the map $t\mapsto \H^2(K_{t_0+t})$ is differentiable, with 
\be \frac{d}{dt}|_{t=t_0}\H^2(K_t)=0.\ee
\qed

\section{Properties for the stabilities}

In this section, we prove two properties :

$1^\circ$ The Almgren sliding stability and the $G$-topological stability of a 2-dimensional minimal cone are independent of the ambient dimension. 

$2^\circ$ The $G$-topological stability implies the $G$-topological sliding stability, which implies Almgren stability, for any 2-dimensional $G$-topological minimal cone (which is automatically Almgren minimal, by Proposition \tb{2.7}). 

These two properties will help to simplify the proof of sliding stability for concrete minimal cones, in the next sections.

\subsection{The Almgren sliding stability and the $G$-topological stability is independent of ambient dimension}

As in Remark \tb{2.8}, we know that if $K\subset \R^m$ is an Almgren (resp. $G$-topological) minimal cone, then it is Almgren (resp. $G$-topologically) minimal of the same dimension in $\R^n$ for all $n\ge m$. That is, the minimality does not depend on the ambient dimension. The following proposition says the same thing for the Almgren sliding stability and the $G$-topological stability.

\begin{rem}We do not know whether the $G$-topological sliding stability is also independent of dimension.
\end{rem}

\begin{pro}Let $K\subset \R^d$ be a 2-dimensional Almgren (resp. $G$-topological) minimal cone. Suppose it is Almgren (resp. $G$-topologically) sliding stable in $\R^d$. Then for all $n\ge d$, $K$ is Almgren (resp. $G$-topologically) sliding stable in $\R^m$.
\end{pro}

\nd Let $K\subset \R^d$ be any 2-dimensional Almgren minimal cone. 

For any $n\ge d$, without loss of generality, suppose that $K$ lies in the $d$-plane $\R^d=\{(x_1,\cdots, x_n): x_{d+1}=\cdots=x_n=0\}$. Let $\cU^n$  denote the $\eta$-convex domain for $K$ in $\R^n$. That is,
\be \cU^n=\{x\in B^n: <x, y>< 1-\eta, \forall y\in K\mbox{ and }<x,a_j><1-2\eta, \forall 1\le j\le m\}\subset \R^n, \ee
where $B^n$ denotes the closed unit ball in $\R^n$.

Also set 
\be A_j^n=\{x\in B^n: <x,a_j>=1-2\eta\mbox{ and }<x,y>\le 1-\eta, \forall y\in K\},\ee
and denote by $\mA^n$ the cone over $A^n:=\cup_{1\le j\le m}A_j^n$ centered at $0$.
Set
\be \Gamma^n_{jl}=\{x\in B^n, <x,y>=1-\eta\mbox{ for some }y\in \gamma_{jl}\}\bs \mA^n,\ee
and
\be S^n_j=\{x\in B^n,<x,y>=1-\eta\mbox{ for some }y\in s^n_j\}.\ee
Set $\G^n=\cup_{(j,l)\in J}\Gamma^n_{jl}$ and $S^n=\cup_{1\le j\le \mu}S^n_j$.

Let $\pi:\R^n\to \R^d$ be the orthogonal projection, which maps a point $x$ to its first $d$ coordinates. Then obviously $\pi(\cU_n)=\cU_d$.

For each $1\le j\le m$, we know that $a_j\in \R^d$ and $K\subset \R^d$. Hence for a point $x\in \R^d$, 
\be\begin{split}
x\in \pi(A_j^n)&\Leftrightarrow
\exists z\in \R^{n-d}\mbox{ such that }(x,z)\in A_j^n\\
&\Leftrightarrow
 \exists z\in \R^{n-d}\mbox{ such that }<(x,z),a_j>=1-2\eta\mbox{ and }<(x,z),y>\le 1-\eta, \forall y\in K\\
 &\Leftrightarrow <x,a_j>=1-2\eta\mbox{, and }<x,y>\le 1-\eta, \forall y\in K\\
 &\Leftrightarrow x\in A_j^d.
\end{split}\ee
Hence $\pi(A_j^n)=(A_j^d)$. Similarly we have
\be \pi(A_j^n)=A_j^d, \forall 1\le j\le m; \pi(\G_j^n)=\G_j^d, \forall (j,l)\in J; \mbox{ and }\pi(S_j^n)=S_j^d,\forall 1\le j\le \mu.\ee

$1^\circ$ We first prove the theorem for the Almgren sliding case.

Let $K\subset \R^d$ be a 2-dimensional Almgren sliding stable minimal cone. Then there exists $\eta<\eta_1(K)$ and $\d<\eta$, such that $K$ is $(\eta,\d)$-sliding stable. We will prove that $K$ is also $(\eta,\d)$-sliding stable in $\R^n$.

First, let $F=\varphi_1(K)$ be a $\d$-sliding deformation in $\cU^n$. The associated family is denoted by $\varphi_t$ as in Definition \tb{2.10}. Let $f:\R^d\to\R^n$ be the inclusion map, i.e. $\pi\circ f=id$. Then the maps $\psi=\pi\circ \varphi_1\circ f$ is a continuous map from $\bar \cU^d$ to $\bar \cU^d$, and $\psi(K)\cap\bar\cU^d=\pi(F\cap\bar\cU^n)$. Moreover it is easy to see that 
\be |\psi(x)-x|<\d\mbox{ for all }x\in K\cap \pa\cU^d.\ee

Now for each $x\in \pa\cU^d\cap K$,  we know that $\varphi_1(x)\in \pa \cU^n$, and $|\varphi_1(x)-x|<\d$. Hence $\varphi_1(x)\in B(K,\d)\cap \pa\cU^d\subset A^d\cup\G^d\cup S^d$. By \tb{(4.6)}, we know that $\pi\circ\varphi_1(x)\in A^d\cup\G^d\cup S^d\subset\pa\cU^d$. As a result, we know that $\psi(K)\subset \pa\cU^d$. We extend $\psi|_K$ to a Lipschitz map $\psi_1$ from $\bar \cU^d$ to $\bar \cU^d$, so that $\psi_1(\pa\cU^d)\subset \pa\cU^d$ and $|\psi_1(x)-x|<\d$ for $x\in \pa\cU^d$.

Next we would like to prove that this $\psi_1$ is a $\d$-sliding Lipschitz deformation of $K$ in $\cU^d$. So let us construct a family $\psi_t,0\le t\le 1$ which satisfies the conditions in Definition \tb{2.10}.  

For $t\in [0,1]$, let $\psi'_t, t\in [0,1]$ be the line homotopy between $\psi_0=id$ and $\psi_1$: $\psi'_t(x)=(1-t)x+t\psi_1(x)$ for $x\in \cU^d$. Then $\psi'_t(\bar\cU^d)\subset \bar\cU^d$, because $\bar \cU^d$ is convex. The problem is that the image of $\pa\cU^d$ under the line homotopy may not stay in $\pa \cU^d$. So we have to push $\psi'_t(\pa \cU^d)$ to $\pa \cU^d$. 

So let $g:B_{\R^d}(\pa\cU^d, (1-t)\d)\to\pa \cU^d$ be a Lipschitz neighborhood retract in $\R^d$. Define, for each $t\in [0,1]$, $h_t:\bar \cU^d\to \bar\cU^d$ as
\be h_t(x)=\left\{\begin{array}{rcl}x&,\ if\ &d(x, \pa\cU^d)\ge2t\d;\\
g(x)&,\ if\ &d(x, \pa\cU^d)\le t\d;\\
d_{t,x}x+(1-d_{t,x})g(x)&,\ if\ &t\d<d(x, \pa\cU^d)<2t\d,\end{array} \right.
\ee
where $d_{t,x}=d(x, \pa\cU^d)/t\d-1\in [0,1]$ for $x$ with $\d<d(x, \pa\cU^d)<2t\d$. 

Now we define, for $x\in \bar\cU$ and $t\in [0,1]$,
\be\psi_t(x)=h_{\frac 12-|t-\frac 12|}\circ \psi'_t.\ee

Then $\psi_t:\bar\cU^d\to \bar\cU^d$ is a continuous family of maps, hence it satisfies \tb{(2.13)}. Note that when $t=0$ or 1, $h_{\frac 12-|t-\frac 12|}=h_0=id$, hence the definition of $\psi_1$ in \tb{(4.9)} coincides with the $\psi_1$ defined previously, and $\psi_0=id$. That is, \tb{(2.12)} is satisfied. By definition of $\psi_1$, \tb{(2.14)} holds.

Now let us look at $\psi_t$ for $0<t<1$. Note that for each $x\in \pa \cU^d$, $\psi_1(x)\in \pa\cU^d$, hence we know that 
\be \begin{split}d(\psi'_t(x), \pa \cU^d)&=d((1-t)x+t\psi_1(x), \pa\cU)\\
&\le \min\{|(1-t)x+t\psi_1(x)-x|, |(1-t)x+t\psi_1(x)-\psi_1(x)|\}\\
&=\min\{t|x-\psi_1(x)|, (1-t)|x-\psi(x)|\}\le\d\min\{t, 1-t\}=\d (\frac 12-|t-\frac 12|).\end{split}\ee
because $|\psi_1(x)-x|<\d$. As a result, $\psi'_t(x)\in B(\pa\cU^d, \d (\frac 12-|t-\frac 12|))$, which means, by definition of $h_t$, that $h_{\frac 12-|t-\frac 12|}(\psi'_t(x))\in \pa \cU^d$. Hence $\psi_t(\pa\cU^d)\subset \pa\cU^d$, which yields \tb{(2.15)}.

Now for \tb{(2.16)}, by definition we have
\be |\psi_t(x)-x|=|(1-t)x+t\psi_1(x)-x|=t|\psi_1(x)-x|\le t\d\le \d.\ee

Now the family $\psi_t,0\le t\le 1$ satisfies \tb{(2.12)-(2.16)}, hence $\psi_1(K)\cap\bar\cU^d=\pi(F\cap\bar\cU^n)$ is a $\d$-sliding deformation of $K$ in $\bar\cU^d$. This holds for any $\d$-sliding deformation of $K$ in $\cU$.

Now let $E\in \oF_\d(K,\bar\cU^n)$. Then by definition, there exists a sequence $\{E_j\}_j$ of $\d$-sliding deformation of $K$ in $\bar\cU^n$, such that $E_j$ converges to $E$. By the above argument, each $\pi(E_j\cap \bar\cU^n)$ is a $\d$-sliding deformation of $K$ in $\bar\cU^d$. Since $\pi$ is Lipschitz, hence $\pi(E_j\cap \bar\cU^n)$ converges to $\pi(E\cap \bar\cU^n)$, and therefore $\pi(E\cap \bar\cU^n)\in\oF_\d(K,\bar\cU^d)$. Since $K$ is $(\eta,\d)$-Almgren sliding stable in $\R^d$, we know that 
\be \H^2(\pi(E\cap \bar\cU^n))\ge\H^2(K\cap\bar\cU^d),\ee
and hence
\be\H^2(K\cap \bar\cU^n)=\H^2(K\cap \bar\cU^d)\le\H^2(\pi(E\cap \bar\cU^n))\le\H^2(E\cap \bar\cU^n).\ee

Note that this holds for all $E\in \oF_\d(K,\bar\cU^n)$. Hence $K$ is $(\eta,\d)$-Almgren sliding stable.

$2^\circ$ Now we prove the theorem for the $G$-topological stability case.

Let $K\subset \R^d$ be a 2-dimensional $G$-topological stable minimal cone. Then there exists $\eta<\eta_1(K)$ and $\d<\eta$, such that $K$ is $(\eta,\d)$-$G$-topological stable. We will prove that $K$ is also $(\eta,\d)$-$G$-topological stable in $\R^n$.

For any $n\ge d$, set 
\be V_\d^n=\{x\in \pa\cU^n:\mbox{ dist}(x,K)\le \d\}.\ee

Take any $(\eta,\d)$-$G$-topological competitor $F$ for $K$ in $\R^n$. This means, by definition, that $F$ satisfies 
\be F\bs\bar\cU^n=K\bs\bar\cU^n,\ee
\be F\cap \pa\cU^n\subset V_\d^n,\ee
and
\be F\cup V_\d^n\mbox{ is a }G\mbox{-topological competitor for }K\mbox{ in }\R^n.\ee

Let us look at the set $\pi(F)$. By \tb{(4.15)}, we know that $F\bs \bar\cU^n=K\bs \bar\cU^n$. Since $K\subset \R^d$, we know that $K\bs \bar\cU^n\subset (\bar \cU^d\times \R^{n-d})^C$, hence $F\bs \bar\cU^n\subset (\bar \cU^d\times \R^{n-d})^C$, and therefore $F\bs \bar\cU^n=F\bs (\bar\cU^d\times \R^{n-d})$. As a result,
\be\pi(F)\bs \bar\cU^d=\pi(F\bs (\bar\cU^d\times \R^{n-d}))=\pi(F\bs\bar\cU^n)=\pi(K\bs\bar\cU^n)=K\bs \bar\cU^d,\ee
which gives \tb{(4.15)} for the set $\pi(F)$ and $n=d$.

Next, since $\d<\eta$, we know that $V_\d^n\subset A^n\cup \G^n\cup S^n$, and hence 
\be \pi(V_\d^n)\subset A^d\cup \G^d\cup S^d\ee
by \tb{(4.6)}. As a result, we have
\be \begin{split}\pi(V_\d^n)&=\pi(\{x\in A^n\cup \G^n\cup S^n:\mbox{dist}(x,K)\le \d\}\\
&\subset \{x\in A^d\cup \G^d\cup S^d:\mbox{dist}(x,K)\le \d\}\\
&=\pi \{x\in A^d\cup \G^d\cup S^d:\mbox{dist}(x,K)\le \d\}\\
&\subset \pi\{x\in A^n\cup \G^n\cup S^n:\mbox{dist}(x,K)\le \d\}=\pi(V_\d^n),\end{split}\ee
which gives
\be\pi(V_\d^n)=\{x\in A^d\cup \G^d\cup S^d:\mbox{dist}(x,K)\le \d\}=V_\d^d.\ee

As a result, $\pi(F)\cap\pa\cU^d=\pi(F\cap (\pa \cU^d\times \R^{n-d}))$. But $F\bs \bar\cU^n=F\bs (\cU^d\times \R^{n-d})$, hence $F\cap \bar\cU^n=F\cap (\cU^d\times \R^{n-d})$. In particular, 
\be F\cap (\pa \cU^d\times \R^{n-d})=F\cap (\cU^d\times \R^{n-d})\cap (\pa \cU^d\times \R^{n-d})=F\cap \bar\cU^n\cap (\pa \cU^d\times \R^{n-d})\subset F\cap \pa\bar\cU^n.\ee
As consequence,
\be \pi(F)\cap\pa\cU^d\subset\pi (F\cap \pa\bar\cU^n)\subset \pi(V_\d^n)=V_\d^d\ee
by \tb{(4.21)}, which gives \tb{(4.16)} for the set $\pi(F)$ and $n=d$. 

Next we prove that $\pi(F)\cup V_\d^d$ is a $G$-topological competitor for $K$ in $\R^d$. Note that $\pi$ is a deformation of $F\cup V_\d^n$ in $2\cU^n$, hence by Proposition \tb{2.7}
$\pi(F\cup V_\d^n)$ is a $G$-topological competitor for $F\cup V_\d^n$ in $\R^n$, and hence by \tb{(4.17)}, 
\be  \pi(F\cup V_\d^n)\mbox{ is a }G\mbox{-topological competitor for }K\mbox{ in }\R^n.\ee

Since $ \pi(F\cup V_\d^n)\subset \R^d$, by Remark \tb{2.8}, we know that 
\be  \pi(F\cup V_\d^n)\mbox{ is a }G\mbox{-topological competitor for }K\mbox{ in }\R^d.\ee

By \tb{(4.23)}, $\pi(F\cup V_\d^n)=\pi(F)\cup \pi(V_\d^n)\subset\pi(F)\cup V_\d^d$, and by \tb{(4.24)}, 
\be \pi(F)\cup V_\d^d\mbox{ is a }G\mbox{-topological competitor for }K\mbox{ in }\R^d.\ee

Combine \tb{(4.18), (4.23) and (4.26)}, we know that $\pi(F)$ is a $(\eta,\d)$-$G$-topological competitor for $K$ in $\R^d$. Since $K$ is $(\eta,\d)$-$G$-topological stable in $\R^d$, we know that 
\be \H^2(\pi(F)\cap\bar\cU^d)\ge \H^2(K\cap\bar\cU^d),\ee
and hence by \tb{(4.18)},
\be \begin{split}\H^2(F\cap \bar\cU^n)&\ge \H^2(\pi(F\cap \bar\cU^n))=\H^2(\pi(F\cap (\bar\cU^d\times\R^{n-d})))\\
&=\H^2(\pi(F)\cap \bar\cU^d)\ge \H^2(K\cap\bar\cU^d).\end{split}\ee

This holds for all $(\eta,\d)$-$G$-topological competitor of $F$ in $\R^n$, hence by definition, $K$ is $(\eta,\d)$-$G$-topological stable in $\R^n$.\qed

\subsection{Topological stable$\Rightarrow$Topological sliding stability$\Rightarrow$Almgren sliding stability}

In this subsection we prove the second property for the stabilities:

\begin{pro}Let $K\subset \R^n$ be a 2-dimensional $G$-topological minimal cone. 

$1^\circ$ If $K$ is $(\eta,\d)$-$G$-topological stable, then it is also $(\eta,\d)$-$G$-topological sliding stable;

$2^\circ$ If $K$ is $(\eta,\d)$-$G$-topological sliding stable, then $K$ is also $(\eta,\d)$-Almgren sliding stable.
\end{pro}

\nd Fix any $\eta$ and $\d$. Let $\cU$ denote $\cU(K,\eta)$. 

$1^\circ$ Let $K\subset \R^n$ be $(\eta,\d)$-$G$-topological stable. We want to prove the $(\eta,\d)$-$G$-topological sliding stability. So let $F$ be a $(\eta,\d)$-$G$-topological sliding competitor for $K$. That is, there exists a 2-dimensional $G$-topological competitor $E$ for $K$ in $\cU$, and a $\d$-sliding deformation $\varphi_t,t\in [0,1]$ for $E$ in $\cU$, such that $F=\varphi_1(E)$. 

We would like to prove that $F$ is a $(\eta,\d)$-$G$-topological competitor for $K$. So let us check the conditions in $1^\circ$ of Definition \tb{2.14}. In fact, \tb{(2.21) and (2.22)} are trivial by definition. So we only have to check \tb{(2.23)}.

We are going to prove that $F':=F\cup V_\d$ is an $G$-topological competitor for $K$ in $4\cU$. We are going to construct a deformation $f$ in $3\cU$, so that $f(E)$ is included in $F'$.

Define $\xi: 3\bar\cU\to 3\bar\cU$:
\be \xi(x)=\left\{\begin{array}{rcl}
x&,\ if\ &x\in \R^n\bs 3 \cU;\\
\frac{3-2(3-|x|)}{|x|}x&,\ if\ &x\in 3 \cU\bs 2 \cU;\\
\varphi_{t_x}(\frac{x}{2-t_x})&,\ if\ &x\in 2 \cU\bs \bar\cU;\\
\varphi_1(x)&,\ if\ &x\in \bar\cU,
\end{array}\right.\ee
where $t_x$ is such that $x\in (2-t_x)\pa\cU$ for $x\in 2 \cU\bs \cU$.

Then $\xi$ is a Lipschitz deformation in $3\cU$, $\xi(E)\bs \bar\cU=E\bs \bar\cU$, $\xi(E\bs 2\cU)\subset E$, $\xi(K\cap \bar\cU)=\varphi_1(E\cap \bar\cU)$, and due to the property \tb{(2.16)} of the family $\varphi_t$, we know that $\xi(E\cap 2\cU\bs \bar\cU)\subset V_\d$.

As a result, $\xi(E)\subset (E\bs \bar\cU)\cup \varphi_1(E\cap \bar\cU)\cup V_\d=F\cup V_\d=F'$.

Note that $\xi(E)$ is a deformation of $E$ in $3\cU$, hence by Proposition \tb{2.7}, $\xi(E)$ is a $G$-topological competitor of dimension 2 for $K$ in $4\cU$. Since $\xi(E)\subset F'$, and $F'\bs 4\cU=\xi(E)\bs 4\cU$, by definition of $G$-topological competitor, $F'$ is also a $G$-topological competitor of dimension $2$ for $E$ in $4\cU$. But $E$ is a $G$-topological competitor for $K$, Hence \tb{(2.23)} holds for $F'$.

As a result, $F$ is a $(\eta,\d)$-$G$-topological competitor for $K$. Since $K$ is $(\eta,\d)$ $G$-topological stable, we know that 
\be\H^2(K\cap\bar\cU)\le \H^2(F\cap \bar\cU).\ee

This holds for all $\d$-sliding $G$-topological competitors $F$ for $K$ in $\cU$, hence $K$ is also $(\eta,\d)$-$G$-topological sliding stable.

$2^\circ$.
Let $K\subset \R^n$ be $(\eta,\d)$-$G$-topological sliding stable. We want to prove that it is $(\eta,\d)$-Almgren sliding stable.  By definition, it is enough to prove that all limits of $(\eta,\d)$-Almgren sliding competitors are $(\eta,\d)$-$G$-topological sliding competitors. But this is easy: if $F=\varphi_1(K)$ be a $\d$-sliding deformation of $K$ in $\cU:=\cU(K,\eta)$, since $K$ is a $G$-topological competitor for itself, we know that $F$ is also a $(\eta,\d)$-$G$-topological sliding competitor. But since Hausdorff limits of $G$-topological competitors in a fixed ball $B$ is still a $G$-topological competitor, we know that all limits of $(\eta,\d)$-Almgren sliding competitors are $(\eta,\d)$-$G$-topological sliding competitors. \qed

\section{Almgren and $G$-topological sliding stability for known 2-dimensional minimal cones contained in $\R^3$}

In this and next section, we are going to prove the Almgren and $G$-topological sliding stabilities for 2-dimensional minima cones in $\R^3$. 

We will see that the proof are getting more and more involved when the structure of the minimal cones are getting more complicated.  To the end, the proof for the codimension 2 cone $Y\times Y$ becomes a separated article (see \cite{stableYXY}).

\subsection{Planes}

\begin{thm}Let $P$ be a $2$-dimensional plane in $\R^n$. Then it is $(\eta,R_1(\eta))$ Almgren and $G$-topological sliding stable for all abelien group $G$, and for all $\eta<1$, where $R_1(\eta)=\sqrt{1-(1-\eta)^2}$.
\end{thm}

\nd Fix any abelien group $G$, and $n\in \N$, and $\eta<1$. We first prove the $G$-topological stability.

Without loss of generality, suppose that $P=\{(x_1,\cdots, x_n)\in \R^n: x_3=x_4\cdots, x_n=0\}$. Let $\pi:\R^n\to P$. Let $\cU$ denote $\cU(P,\eta)$, and let $R_1$ denote $R_1(\eta)$.

We first notice that for any $G$-topological competitor $F$ of $P$ in a ball $B\subset \R^n$, $\pi(F)=P$. In fact, if this is not true, there exists $x\in P$ such that $\pi^{-1}\{x\}\cap F=\emptyset$. Without loss of generality we can suppose that $B=B(0,1)$. Then obviously $x\in B_P(0,1)$. Let $D$ denote the $n-2$-dimensional ball $\{y\in \R^n: \pi(y)=x$ and $||y-x||\le 2\}\subset F^C$, and let $S$ denote the $n-3$-dimensional sphere $S=\{y\in \R^n: \pi(y)=x$ and $||y-x||=2\}\subset F^C$. It is easy to see that $S$ is the boundary of $D$, which does not intersect $F$. Hence $S$ represents a zero element in $H_{n-3}(\R^n\bs F; G)$. Note that $S\subset B(0,1)^C$, $S\cap P=\emptyset$, and $S$ is non zero in $\R^n\bs P$, this contradicts the fact that $F$ is a $G$-topological competitor for $P$.

As a result, if $F$ is a $(\eta,R_1)$-$G$-topological competitor for $P$, then $(F\cup V_{R_1})$ is a $G$-topological competitor for $K$ in $\R^n$, and hence $\pi(F\cup V_{R_1})=P$.

Note that $(F\cup V_{R_1})\bs \bar\cU=P\bs \bar\cU\subset \R^n\bs\pi^{-1}(\bar B_P(0,1-\eta))$, hence $(F\cup V_{R_1})\cap \bar\cU\supset (F\cup V_{R_1})\cap \pi^{-1}(\bar B_P(0,1-\eta))$. Therefore $\pi((F\cup V_{R_1})\cap \bar\cU)\supset \pi^{-1}(\bar B_P(0,1-\eta))$. That is,
\be \pi((F\cap \bar\cU)\cup V_{R_1})\supset\bar B_P(0,1-\eta)=P\cap\bar\cU.\ee
As a result,
\be \begin{split}\H^2(P\cap\bar\cU)&\le \H^2(\pi((F\cap \bar\cU)\cup V_{R_1}))=\H^2(\pi(F\cap \bar\cU)\cup \pi(V_{R_1}))\\&\le\H^2(\pi(F\cap \bar\cU))+\H^2(\pi(V_{R_1}))\le \H^2(F\cap\bar\cU)+\H^2(\pi(V_{R_1})).\end{split}\ee

On the other hand, let $Q=P^\perp$, Then by definition, $V_{R_1}=\{x\in \cU: d(x, P)<R_1\}=\pa B_P(0,1-\eta)\times \bar B_Q(0,R_1)$, hence $\pi (V_{R_1})=\pa B_P(0,1-\eta)$. As a result, 
\be \H^2(\pi (V_{R_1}))=0.\ee

Thus by \tb{(5.2)}, we get
\be \H^2(F\cap\bar\cU)\ge \H^2(P\cap\bar\cU).\ee
This holds for all $(\eta,R_1(\eta))$ $G$-topological competitors for $P$. Hence $P$ is $(\eta,R_1(\eta))$ $G$-topological sliding stable.

By Proposition \tb{4.3}, $P$ is also $(\eta,R_1(\eta))$ $G$-topological sliding and Almgren sliding stable.\qed

\subsection{The $\Y$ set}

In this subsection we prove the sliding stabilities for $\Y$ sets. 

\begin{thm}The $\Y$ sets are $(\eta,R_1(\eta))$-Almgren and $G$-topological sliding stable for all abelien group $G$, and all $\eta<\frac 12$, with $R_1(\eta)=\sqrt{1-(1-\eta)^2}$.
\end{thm}

\nd We first prove the $G$-topological stability. Fix $\eta<\frac 12$. Let $R_1=R_1(\eta)$.

By Proposition \tb{4.2}, it is enough to prove it for the ambient space $\R^3$. So let $Y$ be a 2-dimensional $\Y$ set in $\R^3$ centered at $0$. Without loss of generality, suppose that the spine of $Y$ is the vertical line $\{(x,y,z)\in \R^3:x=y=0\}$, and that the intersection of $Y$ with the horizontal plane $\{z=0\}$ is the union of the three half lines $R_{oa_i},1\le i\le 3$, where $a_1=(1,0,0)$, $a_2=(-\frac 12, \frac{\sqrt 3}{2}, 0)$, and $a_3=-\frac 12, -\frac{\sqrt 3}{2}, 0)$.

Let $\cU$ denote $\cU(Y,\eta)$, let $R_1=R_1(\eta)$. Let $f:\R^3\to \R:f(x,y,z)=z$. Then by definition, $f(\bar\cU)=[-(1-2\eta), 1-2\eta]$. For any set $F\subset \R^3$, and each $t\in \R$, set $F_t=f^{-1}\{t\}\cap F$ the slice of $F$ at level $t$.

We first prove the following lemma:

\begin{lem} If $F$ is a $G$-topological competitor for $Y$ of dimension 2 in a ball $B=B(0,M)$, then for each $t\in [-M,M]$, $F_t\cap B_t$ must connect the three points in $Y_t\cap \pa B=\{a_i^t,1\le i\le 3\}$, i.e.  the three points $a_i^t,1\le i\le 3$ lie in the same connected component of $(F_t\cap B_t)\cup \{a_i^t,1\le i\le 3\}$. \end{lem}

\nd Without loss of generality, suppose that $M=1$. Take any $t\in [-1,1]$. 

Suppose that the three points $a_i^t,1\le i\le 3$ do not belong to the same connected component of $(F_t\cap B_t)\cup \{a_i^t,1\le i\le 3\}$. Suppose for example the connected component $C_1$ of $(F_t\cap B_t)\cup \{a_i^t,1\le i\le 3\}$ contains $a_1^t$ but does not contain $a_2^t$ and $a^t_3$. Then there exists a curve $\gamma:[0,1]\to \bar B_t$ with $\gamma(0),\gamma(1)\in \pa B_t$, which separates $C_1$ and $(F_t\cap B_t)\cup \{a_i^t,1\le i\le 3\}\bs C_1$. That is: 
$ \gamma\subset \bar B_t\bs ((F_t\cap B_t)\cup \{a_i^t,1\le i\le 3\})$, and the sets $C_1$ and $\{a_2^t, a_3^t\}$ belong to different connected components of $\bar B_t\bs \gamma$.  

As  consequence, there exists $t_2, t_3\in [0,1]$, such that $\gamma(t_j)$ belong to the open minor arc of circle $\wideparen{a_1^ta_j^t}$ of $\pa B_t$ between $a_1^t, a_j^t$, $j=2,3$. As a result, $b_j:=\gamma(t_j)$ belong to different connected components of $\R^3_t\bs Y_t$, and hence they belong to different connected components of $\R^3\bs Y$, since $Y=Y_t\times\R$. 

Since $Y$ is a cone, $b_j\not\in Y\Rightarrow$ the segment $[b_j,2b_j]\subset \R^3\bs Y$. Note that $(b_j,2b_j]\subset \R^3\bs \bar B$, and $Y\bs \bar B=F\bs \bar B$, hence $(b_j,2b_j]\subset \R^3\bs F$. Since $b_j\in \bar B_t\bs F_t$, we know that $b_j\in \R^3\bs F$ as well, hence $[b_j,2b_j]\subset \R^3\bs F$. 

Let $\beta$ denote the curve $[2b_2,b_2]\cup \gamma([t_2,t_3])\cup [b_3,2b_3]$. Then $\beta\subset \R^3\bs F$, and it connects $2b_2$ and $2b_3$. Hence the two points $2b_2$ and $2b_3$ belong to the same connected components of $\R^3\bs F$. 

On the other hand, we know that $b_j,j=2,3$ belong to different connected components of $\R^3\bs Y$. Since $[b_j,2b_j]\subset \R^3\bs Y$, $j=2,3$, we know that $2b_j,j=2,3$ belong to different connected components of $\R^3\bs Y$. This contradicts the fact that $F$ is a $G$-topological competitor for $Y$ of codimension 1(which, by Remark 3.2 of \cite{topo}, corresponds to Mumford-Shah competitors, as defined in \cite{DJT} Section 19.)\qed

We continue the proof of Theorem \tb{5.2}. So let $F$ be an $(\eta,R_1)$-$G$-topological competitor of $Y$. By definition, $F':=F\cup V_{R_1}$ is a $G$-topological competitor for $Y$ in a big ball $B=B(0,M)$. Fix any $t\in [-(1-2\eta),1-2\eta]$. By Lemma \tb{5.3}, $F'_t\cap B_t$ must connect $Y_t\cap \pa B$.

Also by definition of $(\eta,\d)$-$G$-topological competitors, we know that $F'_t\bs \bar\cU_t=Y\bs\bar\cU_t$. As a result, $F'_t\cap\bar\cU_t$ connects $Y_t\cap\pa\cU_t$.

Let us look at the planar convex region $\cU_t$. 
%

Denote by $\pi$ the orthogonal projection to the horizontal plane $\{z=0\}$. Then by definition,
\be \pi(\bar \cU_t)=\{q\in \bar B(0,\sqrt{1-t^2}), a_iq\le (1-\eta)^2-t^2, 1\le i\le 3\}.\ee

Note that $\{\sqrt{1-t^2}a_i,1\le i\le 3\}=\pi(Y_t\cap\pa B(0,1))$.

The shape of $\pi(\bar \cU_t)$ is as in the picture: take the disk $B(0, \sqrt{1-t^2})$, then throw away the parts $\{<x,a_i>>(1-\eta)^2-t^2\}$. The boundary of $\pi(\bar \cU_t)$ is a union of three segments $I_i$ of length $R_1$ centered at $[(1-\eta)^2-t^2]a_i$, $1\le i\le 3$, and three arcs of circles $\xi_{ij},1\le i\ne j\le 3$ of $\pa B(0, \sqrt{1-t^2})$ that lie between $I_i$ and $I_j$.

Now we know that $F'_t\cap\bar\cU_t$ connects $Y_t\cap \pa \cU_t$, hence $\pi(F'_t\cap\bar\cU_t)$ connects $\pi(Y_t\cap \pa \cU_t)=\{[(1-\eta)^2-t^2]a_i,1\le i\le 3\}$.
Note that $\pi(F'_t\cap\bar\cU_t)=\pi((V_{R_1}\cap \cU_t)\cup (F_t\cap\bar\cU_t))=\pi(V_{R_1}\cap \cU_t)\cup \pi(F_t\cap\bar\cU_t)$, while $\pi(V_{R_1}\cap \cU_t)=\cup_{1\le i\le 3}I_i$. This means that $[(1-\eta)^2-t^2]a_i,1\le i\le 3$ belong to the same connected componets of $\pi(F_t\cap\bar\cU_t)\cup [\cup_{1\le i\le 3}I_i]$. This implies that $I_i,1\le i\le 3$ belong to the same connected componets of $\pi(F_t\cap\bar\cU_t)\cup [\cup_{1\le i\le 3}I_i]$, because for each $i$, $[(1-\eta)^2-t^2]a_i\in I_i$, which is connected. 

Let $F_0\subset \pi(F_t\cap \bar\cU_t)$ be a subset of $ \pi(F_t\cap \bar\cU_t)$, such that $F_0\cup[\cup_{1\le i\le 3}I_i]$ is connected. Then $F_0$ intersects every $I_i$. For each $i$, fix $q_i\in I_i\cap F_0$. Then $F_0$ is a connected set that contains $q_i,1\le i\le 3$. Let $c\in \Delta_{q_1q_2q_3}$ be the Fermat point of the three points $q_i$. Then we have
\be \H^1(F_0)\ge\sum_{i=1}^3 \H^1([cq_i]).\ee

Denote by $L_i$ the line containing $I_i$, $1\le i\le 3$. Then we have
\be \sum_{i=1}^3 \H^1([cq_i])\ge \sum_{i=1}^3\mbox{dist}(c,I_i)\ge \sum_{i=1}^3\mbox{dist}(c,L_i).\ee

Note that $c$ is contained in the triangle $\Delta_{q_1q_2q_3}$, which is contained in the equilateral triangle $\Delta$ enclosed by the three $L_i,1\le i\le 3$. Since $\Delta$ is equilateral, it is easy to verify that for all $x\in \Delta$, the quantity $\sum_{i=1}^3\mbox{dist}(x,L_i)$ are the same. Hence
\be \sum_{i=1}^3\mbox{dist}(c,L_i)=\sum_{i=1}^3\mbox{dist}(o,L_i)=\sum_{i=1}^3\H^1([0,[(1-\eta)^2-t^2]a_i])=\H^1(\pi(Y_t\cap\bar\cU_t))=\H^1(Y_t\cap\bar\cU_t).\ee

The last equality is because $Y_t$ is parallel to the horizontal plane $\{z=0\}$.

Combine with \tb{(5.6) and (5.7)}, we know that 
\be \H^1(F_0)\ge \H^1(Y_t\cap\bar\cU_t).\ee

Note that $F_0$ is a subset of $\pi(F_t\cap \bar\cU_t)$, hence
\be \H^1(F_0)\le \H^1(\pi(F_t\cap \bar\cU_t))\le \H^1(F_t\cap \bar\cU_t).\ee

As a result, we have
\be \H^1(F_t\cap \bar\cU_t)\ge \H^1(Y_t\cap\bar\cU_t).\ee

Now by coarea formula (cf. \cite{Fe} Theorem 3.2.22), since $f$ is 1-Lipschitz, we have
\be\begin{split}\H^2(F\cap\bar\cU)&\ge \int_{-(1-2\eta)}^{1-2\eta}\H^1(f^{-1}\{t\}\cap (F\cap\bar\cU))=
\int_{-(1-2\eta)}^{1-2\eta}\H^1(F_t\cap \bar\cU_t)\\
&\ge \int_{-(1-2\eta)}^{1-2\eta}\H^1(Y_t\cap \bar\cU_t)=\int_{-(1-2\eta)}^{1-2\eta}\H^1(f^{-1}\{t\}\cap (Y\cap\bar\cU))=\H^2(Y\cap \bar\cU).\end{split}\ee

The last equality is because for almost all $x\in Y$ where the tangent plane $T_xY$ of $Y$ at $x$ exists, $||Df|_{T_xY}||=1$.

Since $F$ is an arbitrary $(\eta, R_1)$-$G$-topological sliding competitor for $Y$, \tb{(5.12)} implies that $Y$ is $(\eta,R_1)$-$G$-topological stable.

By Proposition \tb{4.3}, $Y$ is also $(\eta,R_1(\eta))$ $G$-topological sliding and Almgren sliding stable.\qed

\begin{rem}There are more than one way to prove Theorem \tb{5.2}. The method in the following subsection is more general, which can be used to prove the sliding stabilities for many paired calibrated sets, as $\Y$, $\T$ (see the next section) and $Y\times Y$ (cf. \cite{stableYXY}).
\end{rem}

\subsection{The $\T$ sets}

In this subsection we prove that the 2-dimensional $\T$ sets are Almgren and ($\Z$-)topological sliding stable. Here the Almgren sliding stability will be proved in $\R^n$ for arbitrary $n$, but the topolgical sliding stability is only proved in $\R^3$. The author believes that the topological sliding stability for $\T$ sets are also true in $\R^n$ for any $n\ge 3$. But as stated in the introduction, the main motivation (up to now) of proving the sliding stability for $\T$ sets is to use them to generates new minimal cones by taking the union with other Almgren (resp. topological) sliding stable and Almgren (resp. topological) unique minimal cones. Thus it is essentially enough to know that $\T$ sets are sliding stable in $\R^3$: let $T$ be a $\T$ set, and let $K$ be another  Almgren (resp. topological) unique and  Almgren (resp. topological) sliding stable minimal cone in $\R^d$. Then the almost orthogonal union of $T$ and $K$ in $\R^{3+d}$ is  Almgren (resp. $G$-topological) minimal. Then by Remark \tb{2.8}, the set is also Almgren (resp. topological) minimal in ambient dimension $n+d$ for $n\ge 3$--independent of knowing that $T$ is Almgren (resp. topological) sliding stable in $\R^n$.

Let us first prove the following geometric facts. 

\begin{lem}

Let $P$ be a 2-dimensional plane in $\R^3$. Let $\xi\subset P\cap\partial B(0,1)$ be an arc of angle $\theta<\pi$. Fix a unit vector $v\in P^\perp$, take any $\a\in (0,\frac\pi 2)$, and any $\eta<\frac 12$. Then the following holds: 

Let $e_1$ be the midpoint of $\xi$, and let $e_2\in P\cap\partial B(0,1)$ such that $e_2\perp e_1$. Let $Q_+$ be the 2-plane generated by $e_2$ and $\cos\a e_1+\sin\a v$, and let $Q_-$ be the 2 plane generated by $e_2$ and $\cos\a e_1-\sin\a v$. Let $\pi_\pm$ be the orthogonal projections to $Q_\pm$ respectively.

Let $\t\xi$ be the cone over $\xi$ centered at the origin. Define the 2-dimensional band 
\be \begin{split}&\Xi_\eta:=\xi\times[-R_1 v,R_1 v]\\
=&\{x\in (\t\xi\times P^\perp)\cap \bar B(0,1): <x,z><1-\eta, \forall z\in \xi\mbox{, and }<x,z_0>=1-\eta\mbox{ for some }z_0\in \xi\}.\end{split}\ee

Then for any decomposition $\Xi_\eta=C_+\cup C_-$, with $\H^2(C_+\cap C_-)=0$, we have
\be \H^2(\pi_+(C_+))+\H^2(\pi_-(C_-))=C(\a,\theta)\H^2(\Xi_\eta),\ee
where $C(\a,\theta)$ depends only on $\a$ and the arc length $\theta$ of $\xi$, but not on the decomposition.
\end{lem}

\nd Without loss of generality, we can suppose $P$ is the horizontal plane $\{(x,y,z)\in \R^3:z=0\}$, and $v=(0,0,1)$.

Let $E$ be such a set. Since $\Xi_\eta:=\xi\times[-R_1,R_1]$, for each $q\in \Xi_\eta$, we express it in the coordinate $q=(p,t)$, with $p\in \xi$ and $t\in [-R_1,R_1]$.

For each $p\in \xi$, let $C_{+,p}=C_+\cap (\{p\}\times [-R_1,R_1])$, and $C_{-,p}=C_-\cap (\{p\}\times [-R_1 ,R_1])$. Then we have $\{p\}\times [-R_1,R_1]=C_{+,p}\cup C_{-,p}$. And since $\H^2(C_+\cap C_-)=0$, for almost all $p\in \xi$, $\H^1(C_{+,p}\cap C_{-,p})=0$, and hence
\be \H^1(\{p\}\times [-R_1,R_1])=\H^1(C_{+,p})+\H^1(C_{-,p})\mbox{ for a.e. }p\in \xi.\ee

Since $\Xi_\eta$ is the product of the arc of circle $\xi$ with the interval $[-R_1,R_1]$, at each point $p\in \xi$, the tangent plane $T_{p,t}$ of $\Xi_\eta$ at the point $(p,t)$ is just the plane $p^\perp$ orthogonal to $p$ in $\R^3$. In particular, it does not depend on $t$. The angle between $p^\perp$ and $Q_+$ is the angle between $p$ and $Q_+^\perp=-\sin\a e_1+\cos\a v$. So if $p=\cos\beta_p e_1+\sin\beta_p e_2$ ($\beta\in (-\frac\pi 2,\frac \pi2)$ because $\xi$ is of angle less than $\pi$), then the angle between $p^\perp$ and $Q_+$ is $\arccos|<p,Q_+^\perp>|=\arccos |-\sin\a\cdot \cos\beta_p|=\arccos(\sin\a\cos\beta_p)$.

Similarly, the angle between $p^\perp$ and $Q_-$ is also $\arccos (\sin\a\cos\beta_p)$.

As a result, we know that $|\pi_\pm|_{T_{p,t}}|=\sin\a\cos\beta_p$.

Now let us calculate $\H^2(\pi_+(C_+))$: by the area formula (\cite{Fe} Corollary 3.2.20), we know that
\be \begin{split}\H^2(\pi_+(C_+))&=\int_{C_+}|\pi_+\lf_{T_wC_+}|d\H^2(w)=\int_{p\in \xi}d\H^1(p)\int_{t\in C_{+,p}}|\pi_\pm\lf_{T_{p,t}}|dt\\
&=\int_{p\in \xi}d\H^1(p)\int_{t\in C_{+,p}}\sin\a\cos\beta_pdt=\int_{p\in \xi}d\H^1(p)\sin\a\cos\beta_p\H^1(C_{+,p}).
\end{split}\ee

Similarly we have
\be \H^2(\pi_+(C_-))=\int_{p\in \xi}d\H^1(p)\sin\a\cos\beta_p\H^1(C_{-,p}).\ee

We sum \tb{(5.16) and (5.17)}, and get
\be \begin{split}&\H^2(\pi_+(C_+))+\H^2(\pi_+(C_-))\\
&=\int_{p\in \xi}d\H^1(p)\sin\a\cos\beta_p\H^1(C_{+,p})+\int_{p\in \xi}d\H^1(p)\sin\a\cos\beta_p\H^1(C_{-,p})\\
&=\int_{p\in \xi}d\H^1(p)\sin\a\cos\beta_p[\H^1(C_{+,p})+\H^1(C_{-,p})]\\
&=\int_{p\in \xi}d\H^1(p)\sin\a\cos\beta_p\H^1(\{p\}\times [-R_1,R_1])\\
&=2R_1\sin\a\int_{p\in \xi}d\H^1(p)\cos\beta_p,
\end{split}\ee
the second last equality is after \tb{(5.15)}.

Now since the arclength of $\xi$ is $\theta$, we have
\be \int_{p\in \xi}d\H^1(p)\cos\beta_p=\int_{-\frac\theta2}^{\frac\theta2}\cos tdt=2\sin\frac\theta 2.\ee

As a result, by \tb{(5.18)}, we have
\be \H^2(\pi_+(C_+))+\H^2(\pi_+(C_-))=4R_1\sin\a\sin\frac\theta 2.\ee

And, since $\H^2(\Xi_\eta)=2R_1\times \H^1(\xi)=2R_1\theta$, we know that
\be  \H^2(\pi_+(C_+))+\H^2(\pi_+(C_-))=\frac{2\sin\a\sin\frac\theta 2}{\theta}\H^2(\Xi_\eta).\ee \qed

\begin{lem} For each $\a\in [0,\frac \pi 2)$, the following holds:

Let $Y$ be a 2-dimensional $\Y$ set in $\R^3$, and fix $\eta<\frac 12$. As in the definition of $\cU(K,\eta)$ in Section 2, let $a_1$ and $a_2=-a_1$ be the two singular points of $Y\cap \pa B(0,1)$. Then $Y\cap \pa B(0,1)$ is the union of three half circles $\gamma_j,1\le j\le 3$, all of which join $a_1$ and $a_2$. Let $v_j$ be the mid point of $\gamma_j$, and let $w_j$ be perpendicular to both $v_j$ and $a_1$. Let $P_j$ be the 2-plane generated by $w_j$ and $\sin\a a_1-\cos\a v_j$, $j=1,2,3$. Let $\pi_j$ be the orthogonal projection to $P_j$.  

Then if $E_1,E_2,E_3$ are three essentially disjoint subsets of $A_1$, with $\cup_{j=1}^3 E_j=A_1$, then
\be \sum_{j=1}^3\H^2(\pi_j(E_j))=\cos\a \H^2(A_1).\ee
\end{lem}

\nd The proof of this is almost trivial. Since $A_1$ is part of a plane, it is enough to notice that the angle between $A_1$ and the planes $P_j,1\le j\le 3$ are all $\a$.\qed

Now let us prove the stabilities of $\T$ sets. We first do some simplifications:
  
\begin{defn}Let $U\subset \R^n$ be an open subset such that $\pa U$ is piecewise linear. A closed subset $F$ of $\bar U$ is said to be regular in $\bar U$ if there exists a finite smooth triangulation of $\bar U$ such that $F$ is the support of a 2-simplicial sub-complex of this triangulation. Denote by $\F^R_{\bar U}$ the class of regular subsets in $U$. \end{defn}

\begin{pro}Let $K\subset \R^n$ be a 2-dimensional minimal cone. Fix any $\eta<\eta_1(K)$, and let $\cU$ denote $\cU(K,\eta)$. Then
\be \begin{split}&\inf\{\H^2(F):F\mbox{ is a }(\eta,\d)-G-\mbox{topological sliding competitor for }K\}\\
=&\inf\{\H^2(F):F\mbox{ is a regular }(\eta,\d)-G-\mbox{topological sliding competitor for }K\}.\end{split}\ee
\end{pro}

\nd It is enough to prove that, for any $(\eta,\d)$-$G$-topological sliding competitor $F$ for $K$, and any $\e>0$, there exists a regular $(\eta,\d)$-$G$-topological sliding competitor $F'$ for $K$, such that 
\be \H^2(F')<\H^2(F)+\e.\ee

So take any $(\eta,\d)$-$G$-topological sliding competitor $F$ for $K$, and fix any $\e>0$. Since $F$ is a $G$-topological sliding competitor, there exists a $G$-topological competitor $E$ of $K$ in $\cU$, and a $\d$-sliding deformation $\varphi_t$ in $\bar\cU$, such that $F=\varphi_1(E\cap\cU)$. Since $\cU$ is convex, the segment $[x,\varphi_1(x)]\subset \bar\cU$. For each $x\in \R^n$, let $r_x$ denote the number such that $x\in r_x\pa\cU$ (in other words, $r(x):=r_x$ is the Minkowski functional of the convex set $\cU$). Then without loss of generality, we suppose that $\varphi_t(x)=\frac{(1-t)x+t\varphi_1(x)}{r_{(1-t)x+t\varphi_1(x)}}$, the ''projection'' of the segment $[x,\varphi_1(x)]$ out to $\pa\cU$.

Then it is clear that $|x-\varphi_t(x)|<|x-\varphi_1(x)|<\d$, for all $ t<1$.

So fix a $t_0<1$, such that $\H^2(\varphi([t_0,1]\times(E\cap \pa\cU)))<\frac \e 5$. This is possible, because $E\cap \pa\cU=K\cap\pa\cU$ is 1-rectifiable and hence is of $\H^2$ measure zero. Let $\d_0=\sup_{x\in E\cap \pa\cU}|\varphi_{t_0}(x)-x|$. Since the function $|\varphi_{t_0}(x)-x|<|\varphi_1(x)-x|<\d$ for all $x\in E\cap \pa\cU$, and $E\cap \pa \cU$ is compact, we know that $\d_0<\d$. 

Let $\a>0$ be small, to be decided later. We define $\psi_1:\R^n\to \R^n$, so that 
\be \psi_1(x)=\left\{\begin{array}{rcl}
\varphi_1(\frac{x}{1-\a})&,\ if\ &x\in (1-\a)\bar\cU;\\
\varphi_{(\frac{1-r_x}{\a}) +t_0(1-\frac{1-r_x}{\a})}(\frac{x}{r_x})&,\ if\ &x\in \bar\cU\bs (1-\a)\bar\cU.
\end{array}\right.\ee

Then $\psi_t$ is a $\d_0$ sliding deformation in $\bar\cU$, and we have
\be\begin{split} &\H^2(\psi_1(E)\cap \bar\cU)=\H^2(\psi_1(E\cap\bar\cU))\\
=&\H^2(\psi_1(E\cap(1-\a)\bar\cU))+\H^2(\psi_1(E\cap\bar\cU\bs (1-\a)\bar\cU))\\
=&\H^2(\varphi_1(E)\cap\bar\cU)+\H^2(\varphi([t_0,1]\times(E\cap \pa\cU))\\
< &\H^2(F)+\frac\e5.
\end{split}\ee

Then we apply Theorem 4.3.4 of \cite{Fv}, and get a $n$-dimensional polyhedral complex (which is surely the support of a smooth simplicial sub-complex of a triangluation of $\R^n$) $\cal K$, and a Lipschitz deformation $g:|\K|\to |\K|$, such that 

$1^\circ$ The support of $\K$ covers $(1+\a)\bar\cU$;

$2^\circ$ All the polyhedrons in $\K$ are of diameters less than $\min\{\frac{\d-\d_0}{3}, dist((1+\a)\pa\cU, \cU)\}$;

$3^\circ$ $g|_{(1+\a)\psi_1(E)}$ is a Federer-Fleming projection from $\psi_1(E)$ to $\K^2$. In particular, we have $g(\sigma)\subset \sigma$ for all $\sigma\in \K$;

$4^\circ$ $\H^2(g((1+\a)\psi_1(E)))< \H^2((1+\a)\psi_1(E))+\frac\e5$.

For more details of definitions and notations in the above theorem, see \cite{Fv} for more details (or see \cite{YXY} Section 2 for an explanation in english).

By $2^\circ$ and $3^\circ$, we know that $g\circ\psi(E)$ is a regular subset of $\R^n$, and for each $x\in (1+\a)\psi_1(E)$, $|g(x)-x|\le \frac{\d-\d_0}{3}$, and $|g(x)-x|\le dist((1+\a)\pa\cU, \cU)$. The second control of $|g(x)-x|$ yields that $dist(g(x), \bar\cU)>0$ for $x\in (1+\a)(\psi_1(E)\cap\pa\cU)$. Therefore if we denote by $\pi$ the shortest-distance projection from $\R^n$ to $\bar\cU$, then for each $x\in (1+\a)(\psi_1(E)\cap\pa\cU$, we know that $\pi\circ g(x)\in \pa\cU$. Also, since $g((1+\a)(\psi_1(E)))$ is regular, so is $\pi\circ g((1+\a)(\psi_1(E)))$. Last, for each $x\in (1+\a)(\psi_1(E))$, we know that 
\be\begin{split}|\pi\circ g(x)-x|&\le |\pi\circ g(x)-g(x)|+|g(x)-x|\\
&\le dist(g(x),\bar\cU)+\frac{\d-\d_0}{3}\le |g(x)-\frac{x}{1+\a}|+\frac{\d-\d_0}{3}\\
&\le |g(x)-x|+|x-\frac{x}{1+\a}|+\frac{\d-\d_0}{3}\\
&\le \frac{\d-\d_0}{3}+\frac{\a}{1+\a}+\frac{\d-\d_0}{3}=\frac 23(\d-\d_0)+\frac{\a}{1+\a},\end{split}\ee
and hence for all $x\in \pa\cU$, 
\be\begin{split} &|\pi\circ g((1+\a)\psi_1(x))-x|\\
\le &|\pi\circ g((1+\a)\psi_1(x))-(1+\a)\psi_1(x)|+|(1+\a)\psi_1(x)-(1+\a)x|+|(1+\a)x-x|\\
\le &(\frac 23(\d-\d_0)+\frac{\a}{1+\a})+(1+\a)\d_0+\a.
\end{split}\ee

Now we define $h:\bar\cU\to \bar\cU: h(x)=\pi\circ g((1+\a)\psi_1(x))$. Then it is a $[(\frac 23(\d-\d_0)+\frac{\a}{1+\a})+(1+\a)\d_0+\a]$-sliding deformation of $E$ in $\bar\cU$. Moreover, we know that
\be\begin{split} 
\H^2(h(E))&=\H^2(\pi\circ g((1+\a)\psi_1(E)))\le \H^2(g((1+\a)\psi_1(E)))\\
&< \H^2((1+\a)\psi_1(E))+\frac\e 5<(1+\a)(\H^2(F)+\frac\e 5)+\frac \e5\\
&=\H^2(F)+[\a\H^2(F)+\frac{2+\a}{5}\e].\end{split}\ee

Now we take $\a$, so that
\be [(\frac 23(\d-\d_0)+\frac{\a}{1+\a})+(1+\a)\d_0+\a]<\d,\ee
and 
\be \a\H^2(F)+\frac{2+\a}{5}\e<\e.\ee

Then by \tb{(5.28) and (5.29)}, the map $h$ is a $\d$-sliding deformation for $E$ in $\bar\cU$, and $\H^2(h(E))<\H^2(F)+\e.$

Set $F'=h(E)$. Then $F'$ is a regular $(\eta,\d)$-$G$-topological competitor for $K$ that satisfies \tb{(5.24)}.\qed

\begin{rem}The result can be certainly generalized to general domains $U$ with piecewise linear boundary, and sets $E$ such that $E\cap \pa U$ is regular.
\end{rem}

\begin{pro}Let $K\subset \R^n$ be a 2-dimensional minimal cone. Fix any $\eta$, and any $\d\in (0,\eta)$, let $\cU$ denote $\cU(K,\eta)$. Then 
\be \begin{split}\inf \{\H^2(F): F&\mbox{ is an regular }(\eta,\d)-G-\mbox{ sliding topological competitor of }K\}\\
=\inf \{\H^2(F): F&\mbox{ is an regular }(\eta,\d)-G-\mbox{ sliding  topological competitor of }K\\
&\mbox{ and }F\cap \pa\cU\mbox{ is a Lipschitz deformation of }K\cap \pa\cU\}.\end{split}\ee
\end{pro}

\nd It is enough to prove that, for any regular $(\eta,\d)$-$G$-topological sliding competitor $F$ of $K$, and any $\e>0$, there exists a regular $(\eta,\d)$-$G$-topological sliding competitor $F'$ of $K$, such that 
\be F'\cap \pa\cU\mbox{ is a Lipschitz deformation of }K\cap \pa\cU\mbox{, and }\H^2(F')<\H^2(F)+\e.\ee

So let $F$ be any regular $(\eta,\d)$-$G$-topological sliding competitor of $K$. Then there exists a $G$-topological competitor $E$ of $K$ in $\cU$, and a $\d$-sliding deformation $\varphi_t,t\in [0,1]$ in $\bar \cU$, such that $F=\varphi(E\cap \bar\cU)$.

For any $\gamma>0$ small, let $f:\R^3\to \R^3$, $f(x)=(1-\gamma)x$. Then it is easy to see that $f(E)$ is still a $G$-topological competitor of $K$ in $\cU$. For $x\in \cU$, let $r_x$ be such that $\frac{x}{r_x}\in \pa\cU)$. Let $\psi_t:\bar\cU\to \bar\cU$: 
\be\psi_t(x)=\left\{\begin{array}{rcl}(1-\gamma)\varphi_t(\frac{x}{1-\gamma})&,\ if\ &x\in (1-\gamma)\bar\cU;\\
r_x\varphi_t(\frac{x}{r_x})&,\ if\ &x\in \bar\cU\bs (1-\gamma)\bar\cU.
\end{array}\right.\ee

It is easy to verify that $\psi_1$ is a Lipschitz $\d$-sliding deformation for $f(E)$. In fact, in the region $(1-\gamma)\bar\cU$, the action of $\psi_1$ is just dilate it to $\bar\cU$, apply $\varphi_1$, and then shrink it back to $(1-\gamma)\bar\cU$; while for $\bar\cU\bs (1-\gamma)\bar\cU$, for each $x\in r_x\pa\cU$, we just dilate it to $\pa \cU$, apply $\varphi_1\lf_{\pa\cU}$, and shrink it back to $r_x\pa\cU$. 

As a result, since $f(E)\cap (1-\gamma)\bar\cU=(1-\gamma)E$, we have
\be \psi_1(f(E)\cap (1-\gamma)\bar\cU)=(1-\gamma)\varphi_1(E)=(1-\gamma)F;\ee
 and since $f(E)\cap \bar\cU\bs (1-\gamma)\bar\cU=K\cap \bar\cU\bs (1-\gamma)\bar\cU$, we have
 \be \psi_1(f(E)\cap\bar\cU\bs (1-\gamma)\bar\cU)\mbox{ coincides with the cone over }\varphi_1(K\cap\pa\cU)=\varphi_1(E\cap\pa\cU).\ee
 In particular, $\varphi_1\circ f(E)$ is regular, and \tb{(2.12)-(2.16)} hold for $\psi_t,t\in [0,1]$ and the set $f(E)$. Moreover, $\varphi_1\circ f(E)\cap\pa\cU=\varphi(K\cap\pa\cU)$.

Since $f(E)$ is a $G$-topological competitor of $K$ in $\cU$, by definition, $F':=\varphi_1\circ f(E)$ is a regular $(\eta,\d)$-$G$-topological sliding competitor for $K$, which satisfies that $F'\cap\pa\cU$ is a Lipschitz deformation of $K\cap \pa\cU$.

Let us look at the measure of $F'$. We have 
\be \begin{split}F'&=\psi_1(f(E)\cap\bar\cU)=\psi_1(f(E)\cap (1-\gamma)\bar\cU)\cup\psi_1(f(E)\cap\bar\cU\bs (1-\gamma)\bar\cU)\\
&=(1-\gamma)F\cup \psi_1(f(E)\cap\bar\cU\bs (1-\gamma)\bar\cU).\end{split}\ee
and hence
\be \H^2(F')\le \H^2((1-\gamma)F)+\H^2(\psi_1(f(E)\cap\bar\cU\bs (1-\gamma)\bar\cU)).\ee

We know that $\H^2((1-\gamma)F)=(1-\gamma)^2\H^2(F)\le (1-2\gamma)\H^2(F)$, while for $\psi_1(f(E)\cap\bar\cU\bs (1-\gamma)\bar\cU)$, by \tb{(5.36)} and the definition of $\cU$, a simple calculation give
\be \H^2(\psi_1(f(E)\cap\bar\cU\bs (1-\gamma)\bar\cU))\le \gamma\H^1(\varphi_1(K\cap\pa\cU)).\ee

Combine with \tb{(5.38)}, we get
\be \H^2(F')\le \H^2(F)+\gamma\H^1(\varphi_1(K\cap\pa\cU)).\ee

Let $\gamma=\frac{\e}{\H^1(\varphi_1(K\cap\pa\cU))}$, and we get the conclusion.\qed

Now let us prove the following:

\begin{thm}The $\T$ sets are $(\eta,R_1(\eta))$-Almgren and ($\Z$-)topological sliding stable for all $\eta<\frac 12$ in $\R^3$, with $R_1(\eta)=\sqrt{1-(1-\eta)^2}$.
\end{thm}

\nd We first prove the topological sliding stability. Fix $\eta<\frac 12$. Let $R_1=R_1(\eta)$. 

Let $T$ be a 2-dimensional $\T$ set in $\R^3$ centered at $0$. Denote by $\cU=\cU(T,\eta)$.

Let $a_j,1\le j\le 4$ be the 4 singular points of $T\cap \pa B(0,1)$. Then $\cU\bs T$ is composed of 4 equivalent parts $\Omega'_i, 1\le i\le 4$, and $\cU\bs (A\cup \G)$ (here $S$ does not exist) is composed of 4 equivalent spherical parts $\Omega_i,1\le i\le 4$. The index $i$ is such that the center of $\Omega'_i$ (which is the same as the center of $\Omega_i$) is $-a_i,1\le i\le 4$. Note that $\Omega_i\subset \Omega'_i$, $1\le i\le 4$. 

For each $1\le j\ne l\le 4$, let $\xi_{jl}$ denote the common boundary of $\Omega'_j$ and $\Omega'_l$. Equivalently, if $1\le k\ne i\le 4$ are the other two indices, then $\xi_{jl}=\mc_{ik}\cap \pa\cU$, the Lipschitz curve on $\pa\cU$ that connects $(1-2\eta)a_i$ and $(1-2\eta)a_k$. 

For $1\le j\le 4$, let $\xi_j=\cup_{l\ne j}\xi_{jl}$. Then it is the boundary of $\Omega'_j$.

Take any $(\eta,R_1)$-topological sliding competitor $F$ for $T$. Then there exists a topological competitor $E$ of $T$ in $\cU$, and a sliding deformation $\varphi_1$ in $\bar\cU$, such that $\varphi_1(E\cap \bar\cU)=F$. Note that $E\cap\pa\cU=T\cap\pa\cU=\cup_{1\le j<l\le 4}\xi_{jl}$.

By Propositions \tb{5.8 and 5.10}, we can suppose that $F$ is regular, and that 
\be F\cap\pa\cU=\varphi_1(E\cap\pa\cU)\mbox{, and hence }\H^2(F\cap \pa \cU)=0.\ee

The images $\varphi_1(\xi_j)$, $1\le j\le 4$ are closed Lipschitz curves. Note that $\varphi_1(\pa\cU)\subset \pa\cU$, hence $\varphi_1(\xi_j)$ are still closed curves in $\pa\cU$. 

On $\pa\cU$, if we regard $\xi_j$ and $\Omega'_j$ as simplicial chains with coefficient in $\Z_2$, then the map ${\varphi_1}_*$ maps them to the chains represented by $\varphi_1(\xi_j)$, and the image ${\varphi_1}_*(\Omega'_j)$. Then the support $D_j$ of ${\varphi_1}_*(\Omega'_j)$ is a subset of $\pa\cU$, with $\pa D_j=\varphi_1(\xi_j)$. Since $\cup_{1\le j\le 4}\Omega'_j=\pa\cU$, we have the disjoint union $\pa\cU=\cup_{1\le j\le 4}D_j$. Moreover, since $|\varphi_t(x)-x|<R_1$ for all $t\in [0,1]$, we know that 
\be \Omega_j\subset \Omega'_j\bs B(\xi_j,R_1)\subset D_j\subset \Omega'_j\cup B(\xi_j,R_1)\subset \pa\cU\bs [(\cup_{i\ne j}\Gamma_{ij})\cup A_j].\ee

%

On the other hand, since $E$ is a topological competitor for $T$, we know that $E$ separates the 4 connected components of $\R^3\bs(T\cup \cU)$. Let $C_j$ be the connected component of $\R^3\bs E$ that contains $\Omega'_j$ and with $\xi_j\in \pa C_j$, $1\le j\le 4$. Let $\t C_j=C_j\cap \bar \cU$. Then $\Omega_j\subset \pa \t C_j\subset E\cup \Omega'_j$, and $\pa (\pa \t C_j\bs \Omega'_j)=\pa \Omega'_j=\xi_j$. Therefore, $\xi_j$ represents a zero element in $H_1(E\cap \bar\cU)$. As a result, $\varphi_1(\xi_j)$ represents a zero element in $H_1(F\cap \bar \cU)$. In particular, it represents a zero element in $H_1(F\cap \bar \cU, \Z_2)$. Let $\Sigma_j$ be a $\Z_2$ simplicial chain in $F\cap \bar \cU$ with $\pa\Sigma_j=\xi_j$. To save notations, we still denote by $D_j$ the $\Z_2$-chain associated to $D_j$. Since $\pa D_j=\xi_j$, we know that $\pa(\Sigma_j+D_j)=0$, and hence there exists a simplicial 3-chain in $\bar\cU$ whose boundary is $\Sigma_j+D_j$. Let $O'_j$ the support of this simplicial $\Z_2$ 3-chain, then the boundary of $O'_j$ is a 2 dimensional Lipschitz surface contained in the support of $\Sigma_j+D_j$. 

Let us prove that 
\be D_j\bs \pa O'_j\subset F.\ee  In fact, we know that the support of $\Sigma_j+D_j$ is contained in $O'_j$, $O'_j\subset \bar\cU$, and $D_j\subset \pa\cU$, hence $D_j\cap O'_j$ must lie in the boundary $\pa O'_j$ of $O'_j$. As a result, $D_j\bs \pa O'_j\subset D_j\bs O'_j$. But again the support $ |\Sigma_j+D_j|$ of $\Sigma_j+D_j$ is contained in $O'_j$, hence $D_j\bs O'_j\subset D_j\bs |\Sigma_j+D_j|$, which must be contained in $|\Sigma_j|\cap D_j\subset F\cap D_j$. Hence \tb{(5.43)} holds.

On the other hand, we know that $\H^2(F\cap \pa\cU)=0$, and $D_j\subset \pa\cU$. Hence $\H^2(D_j\cap F)=0$. As a result, by \tb{(5.43)}, we know that $D_j\subset \pa O'_j$, and $\pa O'_j\bs D_j\subset |\Sigma_j|\subset F$. In addition, $O'_j\cap \pa\cU=D_j$.

We would like that the domains $O'_j$ do not intersect each other. So set $O_j=O'_j\bs (\cup_{i\ne j}O'_i, j=1,2,3$. Then the $O_j,1\le j\le 4$ are disjoint. Note that $\pa O_j'\subset \cup_{1\le j\le 4}\pa O_j\subset F\cup\pa \cU$. 

Since the $D_j,1\le j\le 4$ are disjoint, and $O'_j\cap\pa \cU=D_j$, we know that $O_j\cap\pa \cU=D_j$. Set $F_j=\pa O_j\bs D_j$. Then 
\be F_j\subset F\mbox{, and }\pa O_j=F_j\cap D_j.\ee

For each point $x\in D_j$, let $n_j(x)$ be the unit normal vector (which is well defined for almost all $x\in D_j$) pointing towards to the origin. And for each $x\in \pa O_j$, let $v_j(x)$ be the unit normal vector pointing outward to $O_j$. Then if $x\in D_j$, $n_j(x)=-v_j(x)$.

By Stokes, we know that
\be 0=\int_{\pa O_j}<v_j(x),a_j>=\int_{F_j}<v_j(x),a_j>+\int_{D_j}<v_j(x),a_j>,\ee

and hence
\be \int_{F_j}<v_j(x),a_j>=-\int_{D_j}<v_j(x),a_j>=\int_{D_j}<n_j(x),a_j>.\ee

Let $P_i$ be the 2-plane orthogonal to $a_i$, and let $\pi_i$ be the orthogonal projection to $P_i$, then the last term of \tb{(5.46)} is just
\be \int_{D_i}<n_i(x), a_i>d\H^2(x)=\H^2(\pi_i(D_i)).\ee

We sum over $1\le j\le 4$, and get
\be \sum_{j=1}^4\int_{F_j}<v_j(x), a_j>d\H^2(x)=\sum_{j=1}^4\H^2(\pi_i(D_i)).\ee

We first look at the left-hand-side of \tb{(5.48)}. Since $F$ is regular, so is $F'=\cup_{i=1}^4 F_i$. Let $F_{ij}=F_i\cap F_j$, $1\le i\ne j\le 4$. Then for each $x\in F_{ij}$,  $n_i(x)=-n_j(x)$. Since the four $O_j,1\le j\le 4$ are essentially disjoint, for $\H^2$-almost all $x\in F'$, $x$ belongs to at most 2 of the $F_j$. Let $F'_j=F_j\bs\cup_{i\ne j}F_{ij}$. Then we have the essentially disjoint union
\be F'=[\cup_{j=1}^4 F'_j]\cup[\cup_{1\le i<j\le 4}F_{ij}].\ee
 
Thus the left-hand-side of \tb{(5.48)} becomes
\be \begin{split}&\sum_{j=1}^4\int_{F_j}<v_j(x), a_j>d\H^2(x)\\
=&\sum_{1\le j\le 4}\int_{F'_j}<v_j(x), a_j>d\H^2(x)+\sum_{1\le i<j\le 4}\int_{F_{ij}}(<v_i(x), a_i>+<v_j(x), a_j>)d\H^2(x)\\
=&\sum_{1\le j\le 4}\int_{F'_j}<v_j(x), a_j>d\H^2(x)+\sum_{1\le i<j\le 4}\int_{F_{ij}}<v_i(x),a_i-a_j>d\H^2(x)\\
\le&\sum_{1\le j\le 4}\int_{F'_j}||a_j||d\H^2(x)+\sum_{1\le i<j\le 4}\int_{F_{ij}}||a_i-a_j||d\H^2(x)\\
=&\sum_{1\le j\le 4}|a_j|\H^2(F'_j)+\sum_{1\le i<j\le 4}||a_i-a_j||\H^2(F_{ij}).\end{split}\ee

Note that $||a_j||=1,1\le j\le 4$, and $||a_i-a_j||=\frac{2\sqrt 2}{\sqrt 3}$, hence
\be\begin{split} \sum_{j=1}^4\int_{F_j}&<v_j(x), a_j>d\H^2(x)\le \sum_{1\le j\le 4}\H^2(F'_j)+\sum_{1\le i<j\le 4}\frac{2\sqrt 2}{\sqrt 3}\H^2(F_{ij})\\
&\le \frac{2\sqrt 2}{\sqrt 3}[\sum_{1\le j\le 4}\H^2(F'_j)+\sum_{1\le i<j\le 4}\H^2(F_{ij})]
=\frac{2\sqrt 2}{\sqrt 3}\H^2(F'),\end{split}\ee
where the last inequality is because of the disjoint union \tb{(5.49)}.

Let us calculate the right-hand-side of \tb{(5.48)}. Take $i=1$ for example. 

By \tb{(5.42)}, we know that $\Omega_j\subset D_j\subset \pa\cU\bs [(\cup_{i\ne j}\Gamma_{ij})\cup A_j$, that is, $D_j\cap A_j=\emptyset$, and $D_j\cap \Gamma_{ij}=\emptyset$ for all $i\ne j$. But since $\cup_{1\le j\le 4}D_j=\pa\cU$, and $A_j,\Gamma_{ij}$ are subsets of $\pa\cU$, $1\le i\ne j\le 4$, we have the disjoint union
\be A_j=\cup_{i\ne j}(A_j\cap D_i),1\le j\le 4,\ee
and 
\be \Gamma_{ij}=(\Gamma_{ij}\cap D_k)\cup (\Gamma_{ij}\cap D_l)\mbox{ for all permutations }(i,j,k,l)\mbox{ of }(1,2,3,4).\ee

Fix any $1\le j\le 4$, we know that the angle between the plane containing $A_j$ and all the 3 $P_i,i\ne j$ are the same. Hence by Lemma \tb{5.6}, and \tb{(5.52)}, the quantity
\be \sum_{i\ne j}\H^2(\pi_i(D_i\cap A_j))\ee
is a constant that does not depend on $F$. In particular, note that $T$ itself is a topological sliding competitor for $T$, where the $D_j$ for $T$ corresponds to the $\Omega'_j$, $1\le j\le 4$. Hence
\be \sum_{i\ne j}\H^2(\pi_i(D_i\cap A_j))=\sum_{i\ne j}\H^2(\pi_i(\Omega'_i\cap A_j)), 1\le j\le 4.\ee

By the same reason, Lemma \tb{5.5} and \tb{(5.53)} gives that, for all permutations $(i,j,k,l)$ of $(1,2,3,4)$, we have
\be  \H^2(\Gamma_{ij}\cap D_k)+\H^2(\Gamma_{ij}\cap D_l)= \H^2(\Gamma_{ij}\cap \Omega'_k)+\H^2(\Gamma_{ij}\cap \Omega'_l).\ee

Recall that we have the essentially disjoint union
\be \pa\cU=[\cup_{1\le j\le 4}\Omega_j]\cup[\cup_{1\le j\le 4}A_j]\cup[\cup_{1\le i<j\le 4}\Gamma_{ij}],\ee
hence
\be D_i=[\cup_{1\le j\le 4}(D_i\cap\Omega_j)]\cup[\cup_{1\le j\le 4}(D_i\cap A_j)]\cup[\cup_{1\le i<j\le 4}(D_i\cap\Gamma_{ij})].\ee
By \tb{(5.42)}, we have
\be D_i=\Omega_i\cup[\cup_{j\ne i}(D_i\cap A_j)]\cup[\cup_{j,k\ne i,j<k}(D_i\cap\Gamma_{ij})].\ee

Note that by definition, the projections $\pi_i(\Omega_i)$, $\pi_i(A_j), j\ne i$, and $\pi_i(\Gamma_{ij}), j,k\ne i,j<k$ are disjoint, hence we have the disjoint union
\be \pi_i(D_i)=[\pi_i(\Omega_i)]\cup[\cup_{j\ne i}\pi_i(D_i\cap A_j)]\cup[\cup_{j,k\ne i,j<k}\pi_i(D_i\cap\Gamma_{ij})],\ee
which gives
\be \H^2(\pi_i(D_i))=\H^2[\pi_i(\Omega_i)]+[\sum_{i\ne j}\H^2(\pi_i(D_i\cap A_j))]+[\sum_{j\ne i,k\ne i,j<k}\H^2(\pi_i(D_i\cap\Gamma_{jk}))].\ee
We sum over $1\le i\le 4$, and get
\be \begin{split}&\sum_{1\le i\le 4}\H^2(\pi_i(D_i))\\
=&\sum_{1\le i\le 4}\{[\H^2(\pi_i(\Omega_i))+[\sum_{i\ne j}\H^2(\pi_i(D_i\cap A_j))]+[\sum_{j\ne i,k\ne i,j<k}\H^2(\pi_i(D_i\cap\Gamma_{jk}))]\}\\
=&[\sum_{1\le i\le 4}[\H^2(\pi_i((\Omega_i))]+[\sum_{1\le i<j\le 4}\H^2(\pi_i((D_i\cap A_j))]+[\sum_{j\ne i,k\ne i,j<k}\H^2(\pi_i((D_i\cap\Gamma_{jk}))]\\
=&[\sum_{1\le i\le 4}[\H^2(\pi_i(\Omega_i))]+[\sum_{1\le j\le 4}(\sum_{i\ne j}\H^2(\pi_i(D_i\cap A_j)))]+[\sum_{j<k}(\sum_{i\ne j,k}\H^2(\pi_i(\Gamma_{jk}\cap D_i))].
\end{split}\ee

The same argument holds for $D_i=\Omega_i$, because $T$ is a competitor, and the $D_i$ correspond to $\Omega'_i$, $1\le i\le 4$. Hence we have
\be \begin{split}&\sum_{1\le i\le 4}\H^2(\pi_i(\Omega'_i))\\
=&[\sum_{1\le i\le 4}[\H^2(\pi_i(\Omega_i))]+[\sum_{1\le j\le 4}(\sum_{i\ne j}\H^2(\pi_i(\Omega'_i\cap A_j)))]+[\sum_{j<k}(\sum_{i\ne j,k}\H^2(\pi_i(\Gamma_{jk}\cap \Omega'_i))].\end{split}\ee

By \tb{(5.55) and (5.56)}, the right-hand-sides of \tb{(5.62) and (5.63)} are equal. Hence we have
\be \sum_{1\le i\le 4}\H^2(\pi_i(D_i))=\sum_{1\le i\le 4}\H^2(\pi_i(\Omega'_i)).\ee

Combine \tb{(5.48), (5.51) and (5.64)}, we have
\be \H^2(F')\ge\frac{\sqrt 3}{2\sqrt 2}[\sum_{1\le i\le 4}\H^2(\pi_i(\Omega'_i))].\ee

On the other hand, either by chasing the condition of equality for the inequalities of \tb{(5.50) and (5.51)}, or by a direct calculus, it is easy to see that
\be \H^2(T\cap \bar\cU)=\frac{\sqrt 3}{2\sqrt 2}[\sum_{1\le i\le 4}\H^2(\pi_i(\Omega'_i))].\ee

Hence we have
\be \H^2(T\cap \bar\cU)\le \H^2(F')\le \H^2(F).\ee
This holds for an arbitrary $(\eta,R_1)$ topological sliding competitor $F$ for $T$, hence $T$ is $(\eta,R_1)$ topological sliding stable. 

The Almgren sliding stability in $\R^3$ follows directly from Proposition \tb{4.3}. This finishes the proof of Theorem \tb{5.11}.\qed
 
 \begin{cor}The $\T$ sets are $(\eta,R_1(\eta))$-Almgren sliding stable for all $\eta<\frac 12$ in $\R^n$ for all $n\ge 3$, with $R_1(\eta)=\sqrt{1-(1-\eta)^2}$.
 \end{cor}
 
 \nd This follows directly from Theorem \tb{5.11} and  Proposition \tb{4.2}.\qed
 
 \renewcommand\refname{References}
\bibliographystyle{plain}
\bibliography{reference}

\end{document}